\documentclass[a4paper,12pt]{article}
\title{
Zalcman functions
and
similarity between the Mandelbrot set, Julia sets, 
and the tricorn
}

\author{
Tomoki Kawahira\\
Tokyo Institute of Technology
\thanks{Partially supported by JSPS}
}%author
\usepackage{enumerate}
\usepackage{graphicx}
\usepackage{amsmath}
\usepackage{amscd}
\usepackage{amssymb}
\usepackage{amsfonts}
\usepackage{theorem}

\newtheorem{thm}{Theorem}%[section]
\newtheorem{prop}[thm]{Proposition}
\newtheorem{lem}[thm]{Lemma}

\theorembodyfont{\rmfamily}
\newtheorem{pf}{Proof.}

\newcommand{\thmref}[1]{Theorem \ref{#1}}
\newcommand{\propref}[1]{Proposition \ref{#1}}

\newcommand{\lemref}[1]{Lemma \ref{#1}}

\newcommand{\parag}[1]{
\medskip
\noindent {\bf #1}
}%---pagagraph 

\newcommand{\C}{\ensuremath{\mathbb{C}}}

\newcommand{\Cstar}{\ensuremath{\mathbb{C}^\ast}}
\newcommand{\Chat}{\ensuremath{\widehat{\mathbb{C}}}}
\newcommand{\R}{\ensuremath{\mathbb{R}}}
\newcommand{\D}{\ensuremath{\mathbb{D}}}

\newcommand{\Z}{\ensuremath{\mathbb{Z}}}
\newcommand{\N}{\ensuremath{\mathbb{N}}}
\newcommand{\M}{\ensuremath{\mathbb{M}}}

\newcommand{\T}{\ensuremath{\mathbb{T}}}
\newcommand{\paren}[1]{{\left( #1 \right)}}
\newcommand{\braces}[1]{{\left\{ #1 \right\}}}
\newcommand{\gauss}[1]{{\left [ #1 \right ]}}
\newcommand{\id}{\mathrm{id}}
\newcommand{\lam}{{\lambda}}

\newcommand{\del}{{\partial}}

\newcommand{\zbar}{\Bar{z}}

\newcommand{\e}{\epsilon}
\newcommand{\cc}{\circ}
\newcommand{\dz}{\, dz}
\newcommand{\cF}{{\mathcal{F}}}

\newcommand{\cJ}{{\mathcal{J}}}

\newcommand{\cM}{{\mathcal{M}}}

\newcommand{\cP}{{\mathcal{P}}}
\newcommand{\cQ}{{\mathcal{Q}}}

\newcommand{\cU}{{\mathcal{U}}}

\newcommand{\cZ}{{\mathcal{Z}}}

\newcommand{\ee}{~=~}
\newcommand{\dee}{~:=~}
\newcommand{\st}{\,:\,}

\newcommand{\QED}{\hfill $\blacksquare$}
\renewcommand{\Bar}{\overline}

\newcommand{\Aff}{\mathrm{Aff}}
\newcommand{\cbar}{\overline{c}}

\newcommand{\cUbar}{\overline{\cU}}

\begin{document}

\maketitle

	\begin{center}
{\it Dedicated to Lawrence Zalcman
 on the occasion of his 75th birthday}
\end{center}

\begin{abstract}%%abstract text
We present a simple proof of Tan's theorem on 
asymptotic similarity between the Mandelbrot set and Julia sets at Misiurewicz parameters. Then we give a new perspective on this phenomenon in terms of Zalcman functions, 
that is, entire functions generated by applying 
Zalcman's lemma to complex dynamics. 
We also show asymptotic similarity between the tricorn 
and Julia sets at Misiurewicz parameters,
which is an antiholomorphic counterpart of Tan's theorem.
\end{abstract}

%-------------------------------------------------

\section{Similarity between $\M$ and $J$}\label{sec_MJ}
The aim of this paper is to give a new perspective on a well-known 
similarity between the Mandelbrot set and Julia sets 
(Tan's theorem) in terms of 
Zalcman's rescaling principle in non-normal families 
of meromorphic functions. 
We start with a simplified proof of Tan's theorem \cite{TL}
following \cite{Ka}, which motivates the whole idea of this paper.

\parag{The Mandelbrot set and the Julia sets.}
Let us consider the quadratic family 
$$
\braces{f_c(z) = z^2 + c \st c \in \C}.
$$ 
The \text{\it Mandelbrot set} $\M$ is the set of $c \in \C$ 
such that the sequence $\braces{f_c^n(c)}_{n \in \N}$ is bounded. 
For each $c \in \C$, the \text{\it filled Julia set} $K_c$ 
is the set of $z \in \C$ such that the sequence 
$\braces{f_c^n(z)}_{n \in \N}$ is bounded.
One can easily check that
\begin{itemize}
\item
$c \notin \M$ if and only if 
$|f_c^n(c)| > 2$ for some $n \in \N$; and
\item
for each $c \in \M$, $z \notin K_c$ if and only if 
$|f_c^n(z)| > 2$ for some $n \in \N$.
\end{itemize}
The {\it Julia set} $J_c$ is the boundary of $K_c$.
Note that all $\M, K_c$, and $J_c$ 
are compact, and also non-empty because 
we can always solve the equations $f_c^n(c)=c$ and $f_c^n(z)=z$.

%%%%%%----------table of pictures
\fboxsep=0pt
\fboxrule=.5pt
\begin{figure}[htbp]
\begin{center} 
\fbox{\includegraphics[width=.18\textwidth, bb = 0 0 600 600]{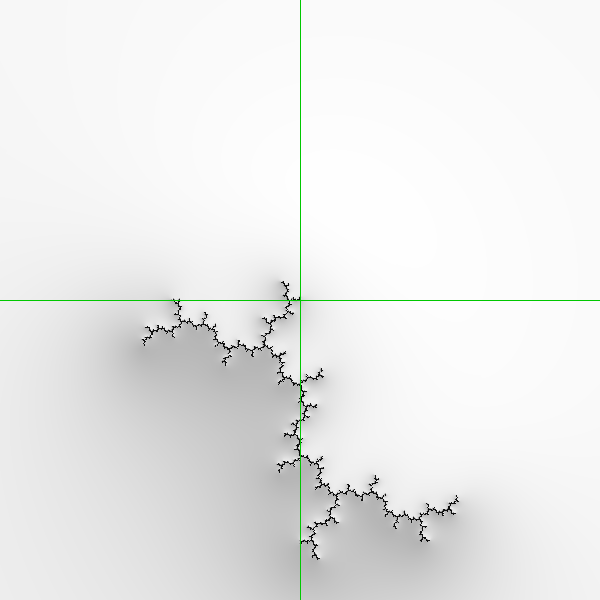}}
\fbox{\includegraphics[width=.18\textwidth, bb = 0 0 600 600]{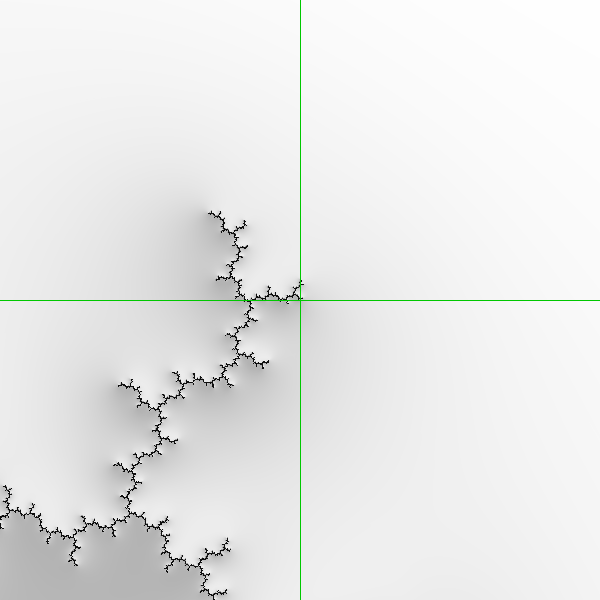}}
\fbox{\includegraphics[width=.18\textwidth, bb = 0 0 600 600]{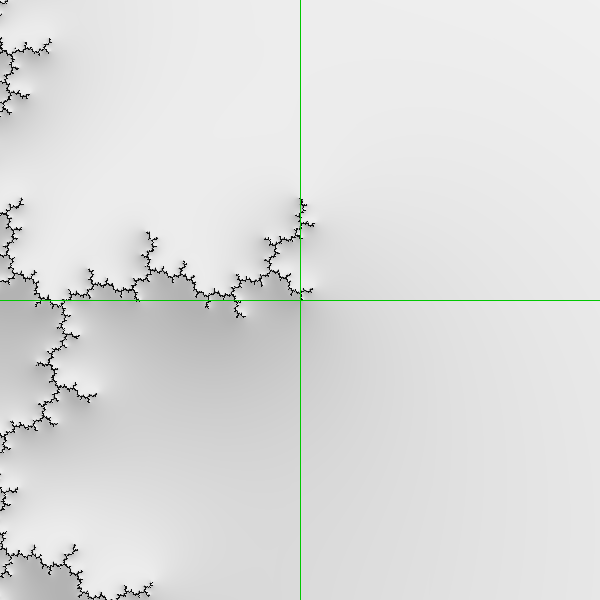}}
\fbox{\includegraphics[width=.18\textwidth, bb = 0 0 600 600]{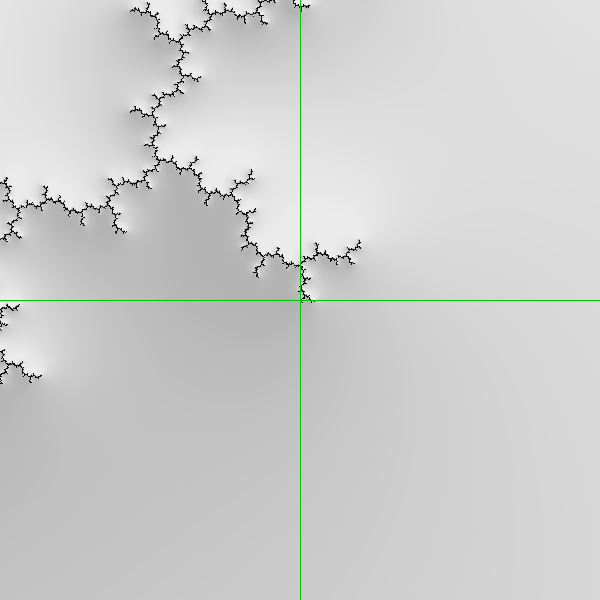}}
\fbox{\includegraphics[width=.18\textwidth, bb = 0 0 600 600]{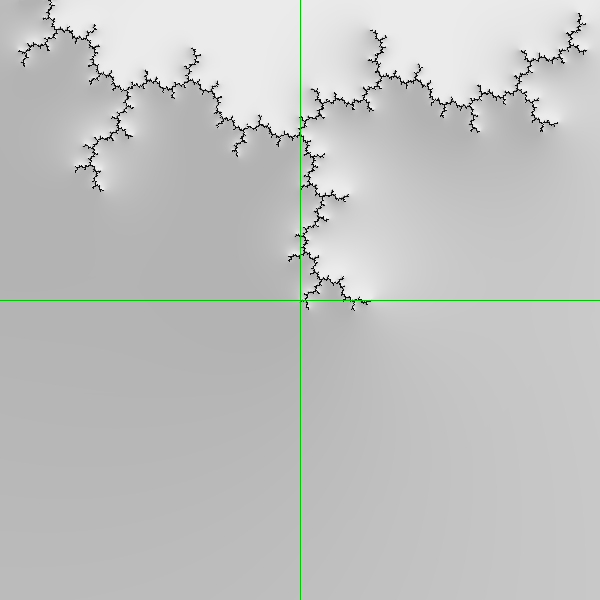}}
\\[.5em]
\fbox{\includegraphics[width=.18\textwidth, bb = 0 0 600 600]{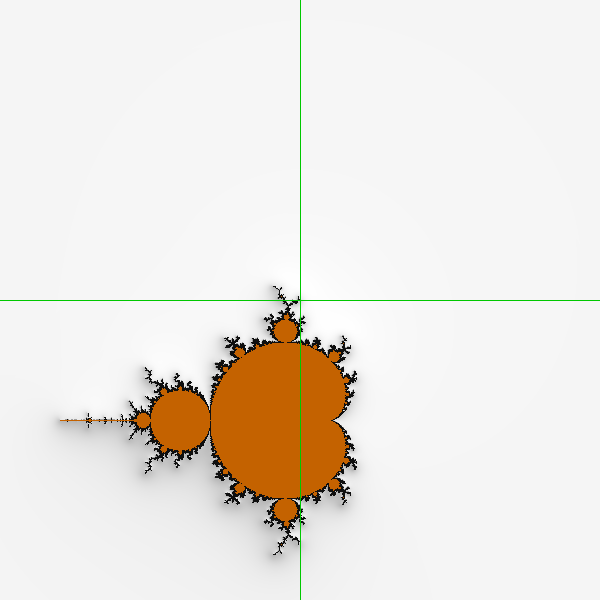}}
\fbox{\includegraphics[width=.18\textwidth, bb = 0 0 600 600]{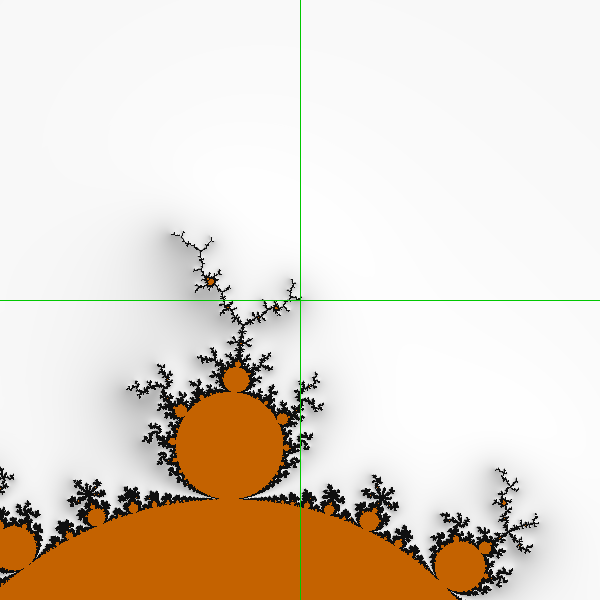}}
\fbox{\includegraphics[width=.18\textwidth, bb = 0 0 600 600]{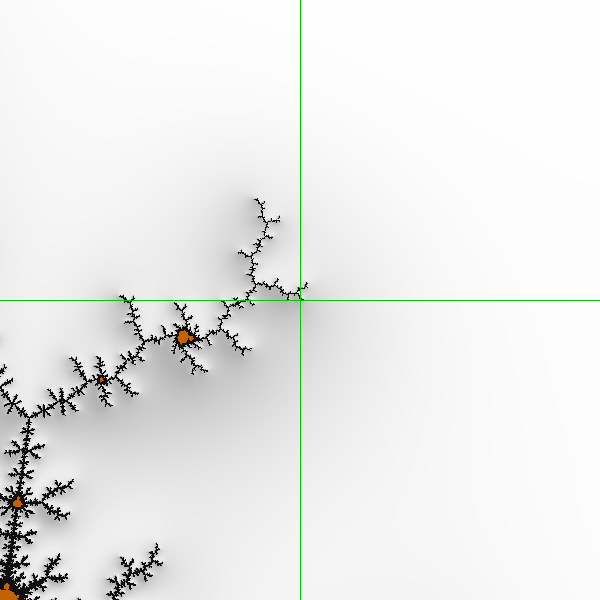}}
\fbox{\includegraphics[width=.18\textwidth, bb = 0 0 600 600]{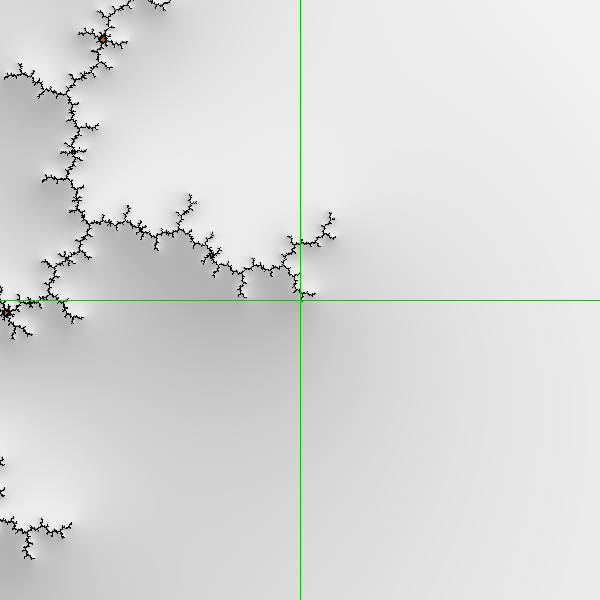}}
\fbox{\includegraphics[width=.18\textwidth, bb = 0 0 600 600]{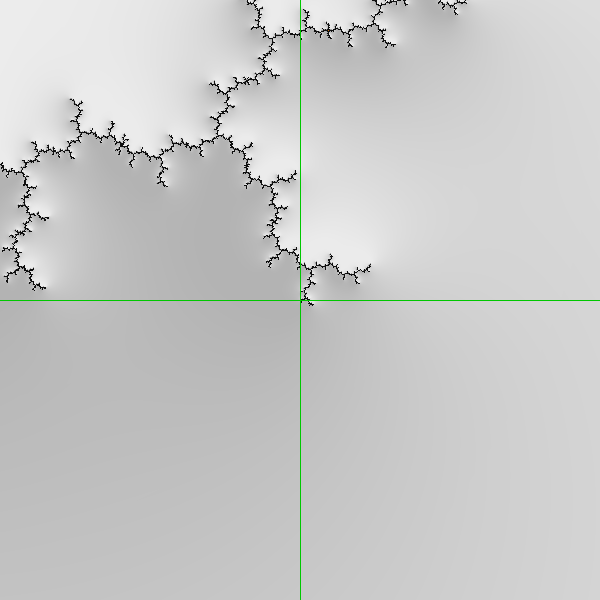}}
\\[.5em]
(JM1)\\[.5em]%:Center: $0.0+1.0 i$, Square width: from 5.0 to 0.01. \\[.5em]
\fbox{\includegraphics[width=.18\textwidth, bb = 0 0 600 600]{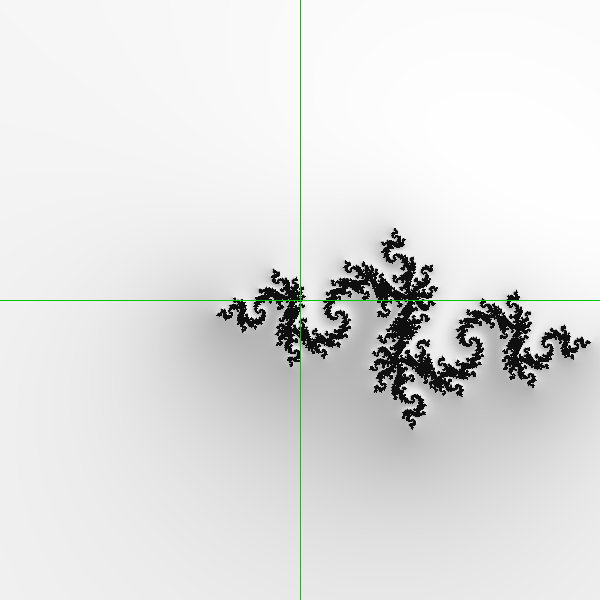}}
\fbox{\includegraphics[width=.18\textwidth, bb = 0 0 600 600]{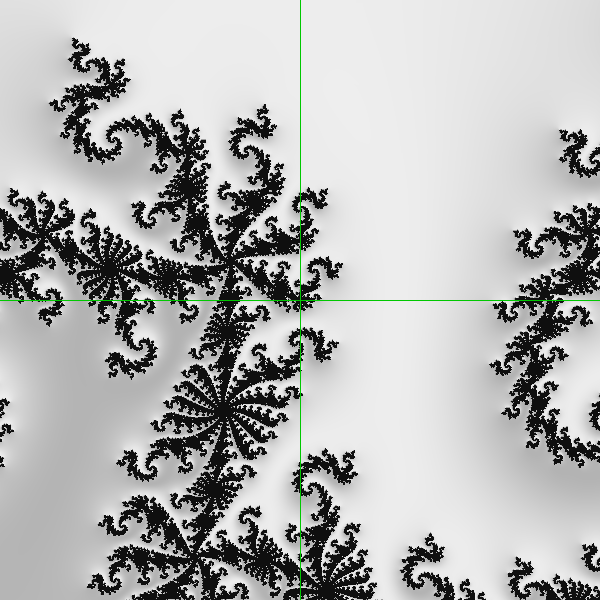}}
\fbox{\includegraphics[width=.18\textwidth, bb = 0 0 600 600]{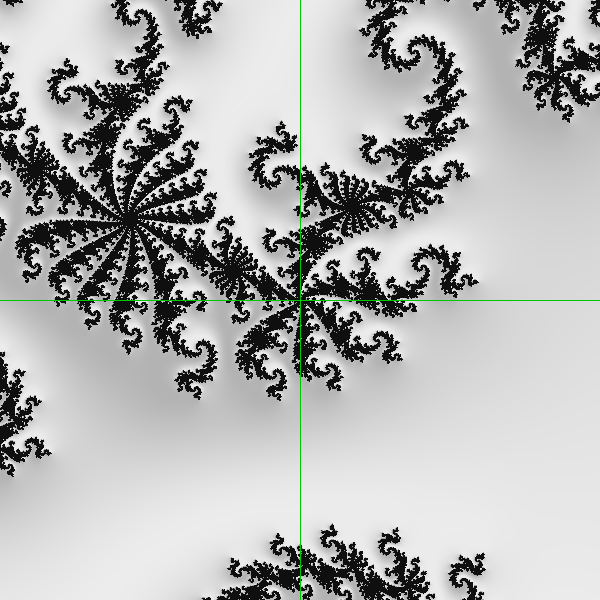}}
\fbox{\includegraphics[width=.18\textwidth, bb = 0 0 600 600]{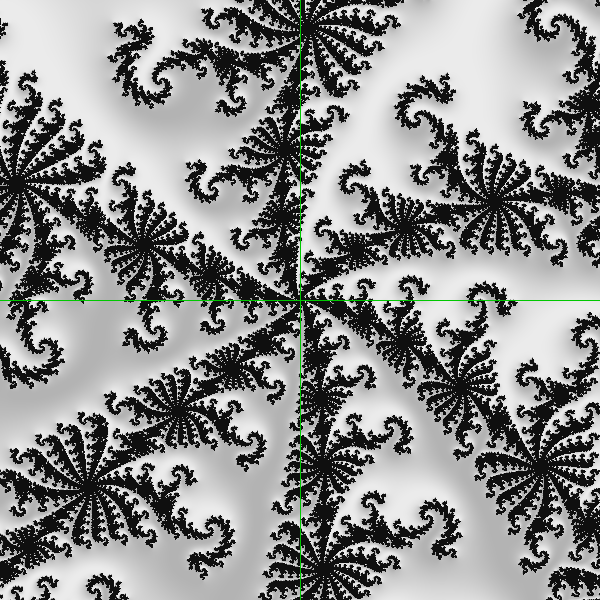}}
\fbox{\includegraphics[width=.18\textwidth, bb = 0 0 600 600]{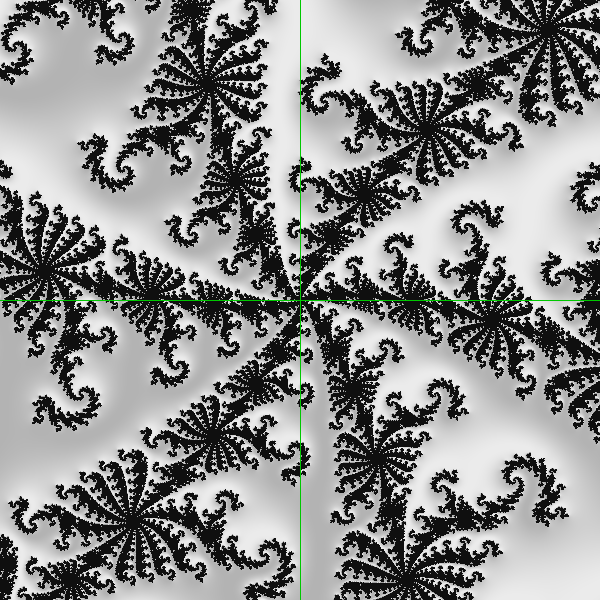}}
\\[.5em]
\fbox{\includegraphics[width=.18\textwidth, bb = 0 0 600 600]{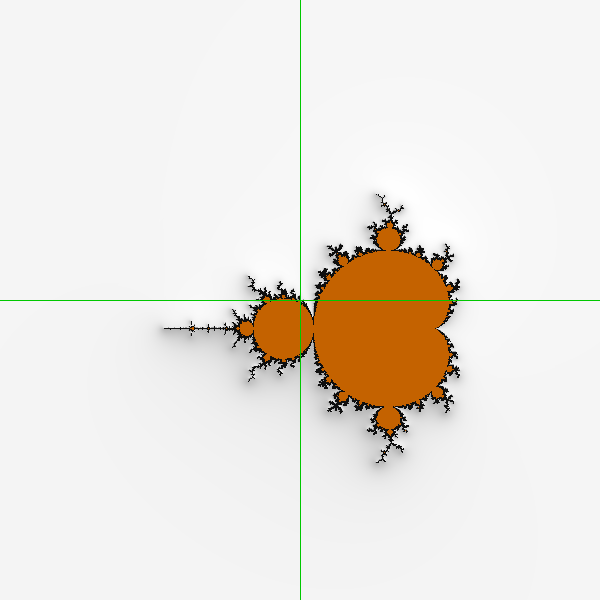}}
\fbox{\includegraphics[width=.18\textwidth, bb = 0 0 600 600]{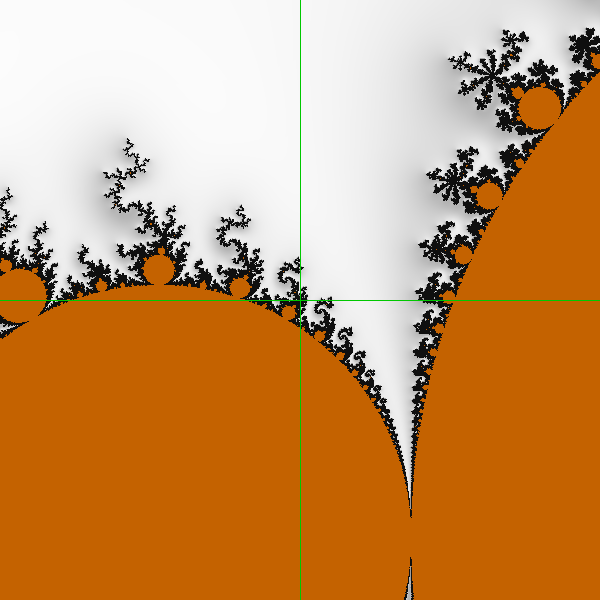}}
\fbox{\includegraphics[width=.18\textwidth, bb = 0 0 600 600]{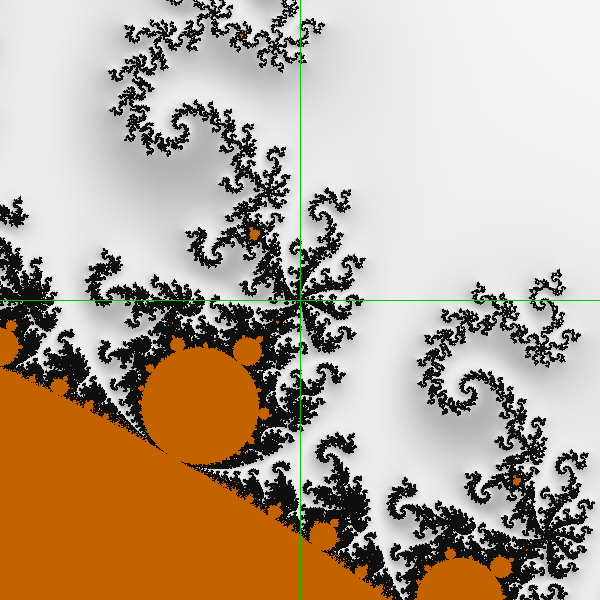}}
\fbox{\includegraphics[width=.18\textwidth, bb = 0 0 600 600]{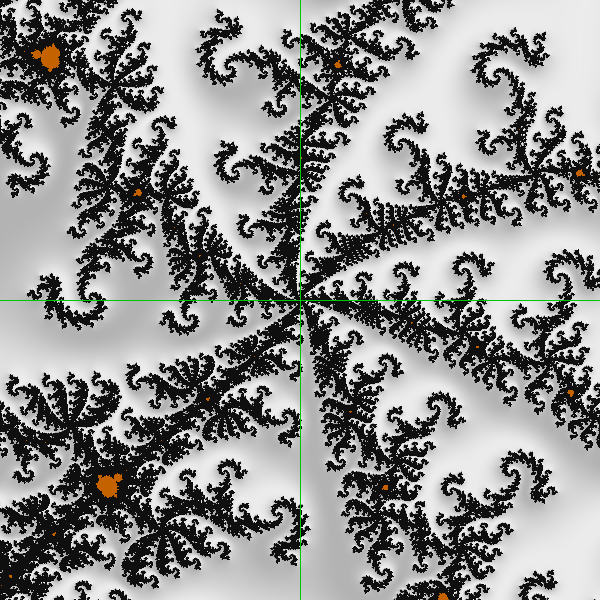}}
\fbox{\includegraphics[width=.18\textwidth, bb = 0 0 600 600]{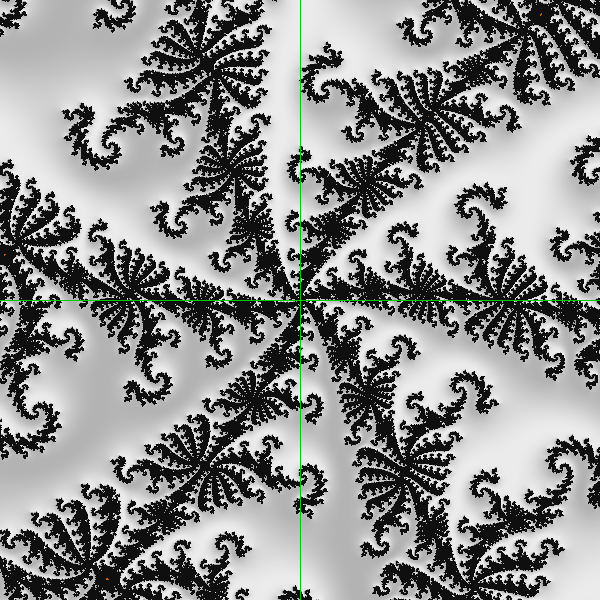}}
\\[.5em]
(JM2)\\[.5em]%:Center: $-0.8597644816892409+0.23487923150145784 i$, 
\fbox{\includegraphics[width=.18\textwidth, bb = 0 0 600 600]{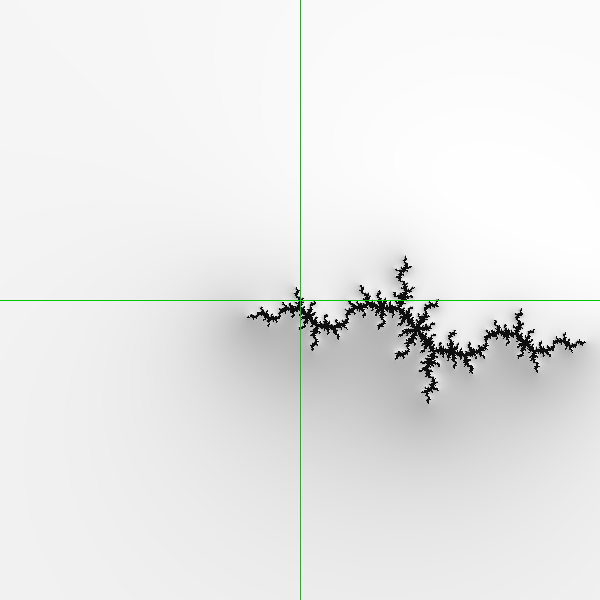}}
\fbox{\includegraphics[width=.18\textwidth, bb = 0 0 600 600]{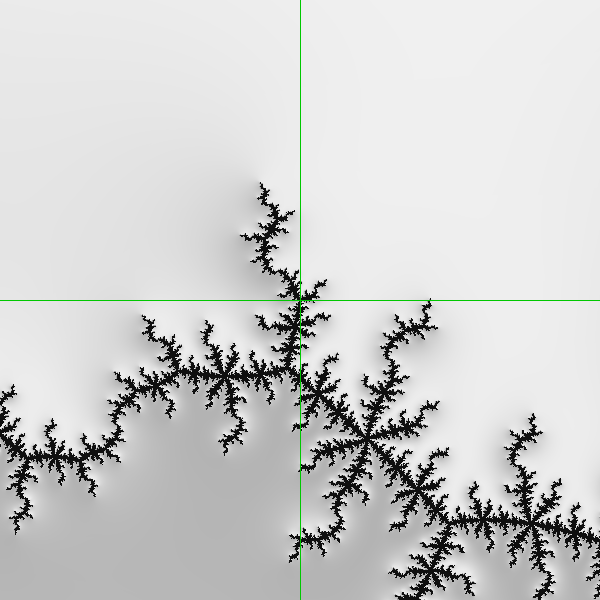}}
\fbox{\includegraphics[width=.18\textwidth, bb = 0 0 600 600]{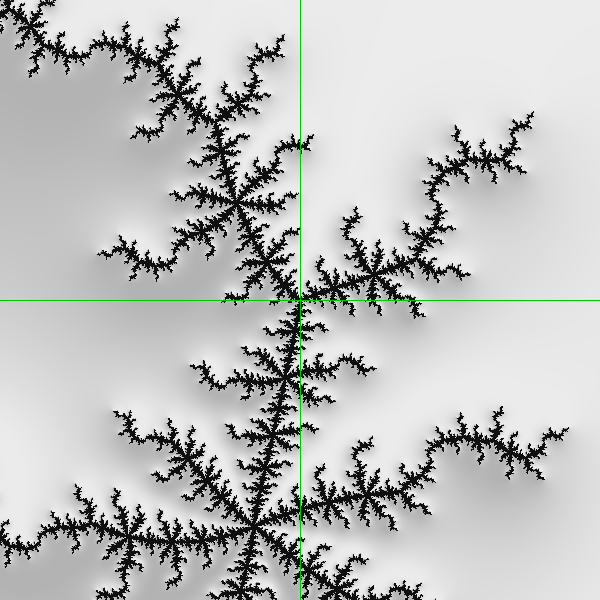}}
\fbox{\includegraphics[width=.18\textwidth, bb = 0 0 600 600]{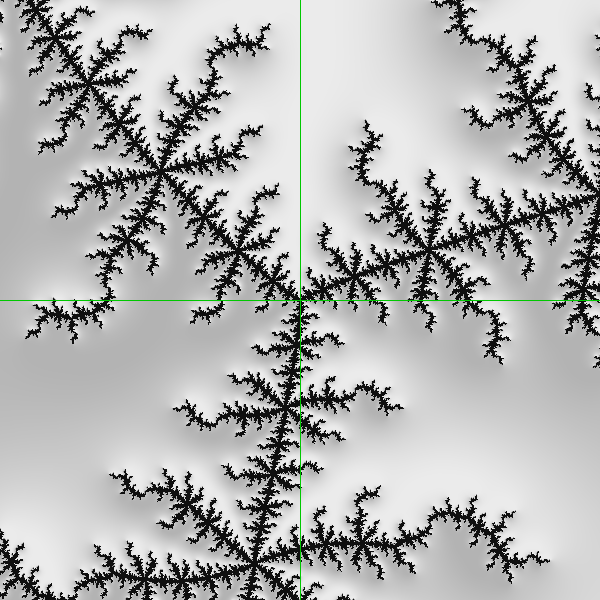}}
\fbox{\includegraphics[width=.18\textwidth, bb = 0 0 600 600]{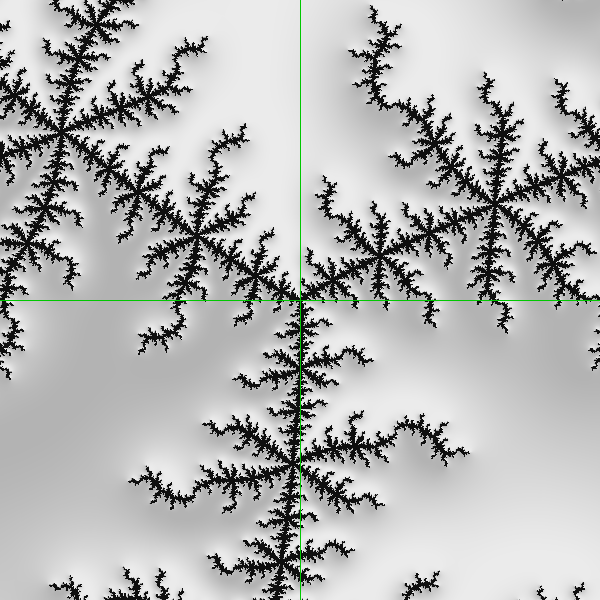}}
\\[.5em]
\fbox{\includegraphics[width=.18\textwidth, bb = 0 0 600 600]{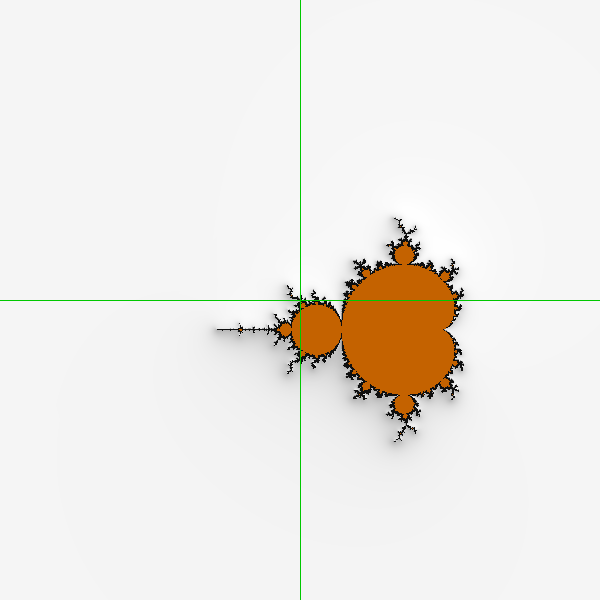}}
\fbox{\includegraphics[width=.18\textwidth, bb = 0 0 600 600]{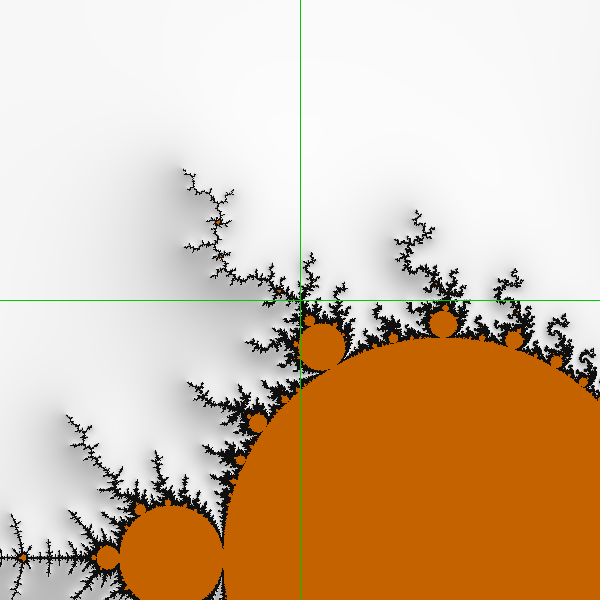}}
\fbox{\includegraphics[width=.18\textwidth, bb = 0 0 600 600]{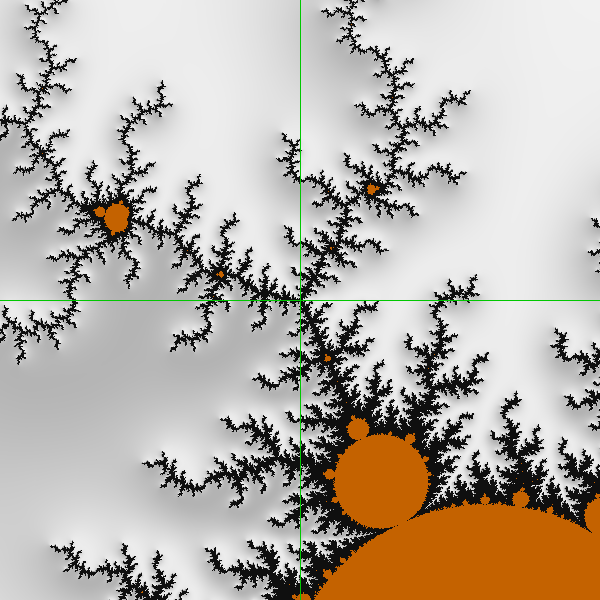}}
\fbox{\includegraphics[width=.18\textwidth, bb = 0 0 600 600]{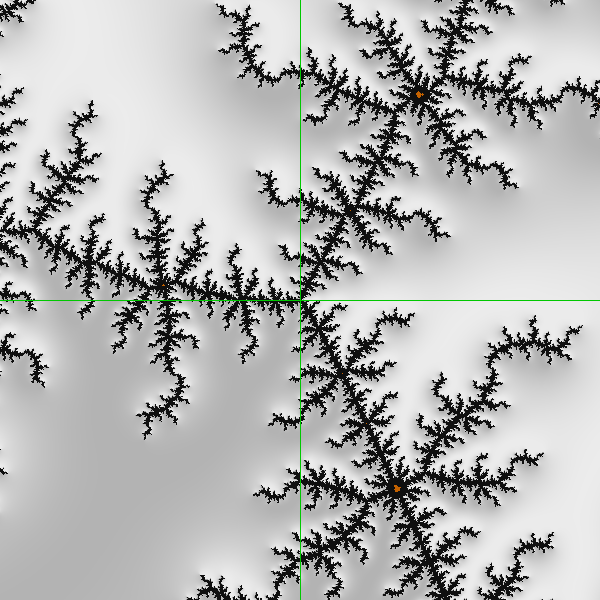}}
\fbox{\includegraphics[width=.18\textwidth, bb = 0 0 600 600]{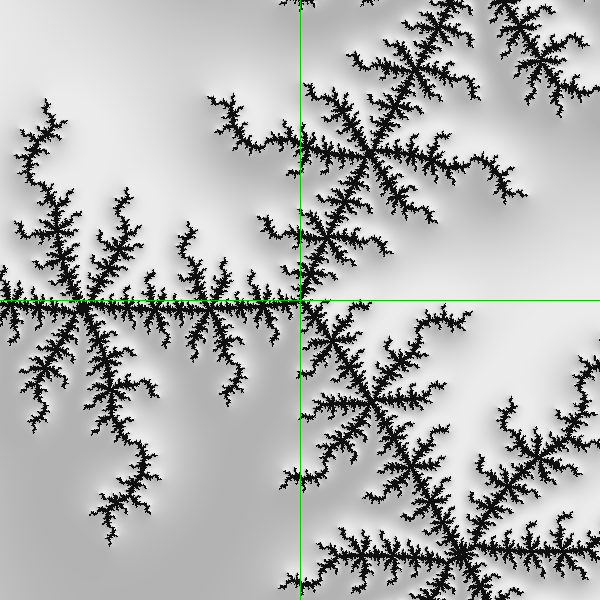}}
\\[.5em]
(JM3)%:Center: $-0.8597644816892409+0.23487923150145784 i$, 
%Square width: from 5.0 to 0.001. 
\\[.5em]
\end{center}
\caption{
(JM1) Center: $0.0+1.0 i$, square width: from 6.0 to 0.01.
(JM2) Center: $-0.8597644816892409+0.23487923150145784 i$, 
square width: from 5.0 to 0.001. 
(JM3) Center: $-1.162341599884035+0.2923689338965703 i$, 
square width: from 6.0 to 0.001.
}
\label{fig_JM}
\end{figure}
Tan showed in \cite{TL} (originally in a chapter of \cite{DH1})
that when $c_0 \in \M$ is a Misiurewicz parameter 
(to be defined below), the ``shapes" of $\M$ and the Julia set $J_{c_0}$ 
are asymptotically similar at the same point $c_0$. 
For example, (JM1) of Figure 1 shows $\M$ and $J_{c_0}$ in squares 
centered at $c_0=i$, whose widths range from $6.0$ to $0.01$.
We will prove this by finding an entire function that bridges the dynamical plane 
and the parameter plane (\lemref{lem_key}).

\parag{Misiurewicz parameters and a key lemma.}
Following the terminology of \cite{DH1} and \cite{TL}, 
we say $c_0 \in \M$ is a \text{\it Misiurewicz parameter} 
if the forward orbit of $c_0$ by $f_{c_0}$ eventually lands on 
a repelling periodic point. 
More precisely, there exist minimal $l \ge 1$ and $p \ge 1$
such that $a_0: = f_{c_0}^l(c_0)$ satisfies 
$a_0 = f_{c_0}^{p}(a_0)$ and $|(f_{c_0}^p)'(a_0)| >1$.
By the implicit function theorem, 
 we can show that the repelling periodic point $a_0$ is stable: 
that is, there exists a neighborhood $V$
of $c_0$ and a holomorphic map $a:c \mapsto a(c)$ on $V$
such that $a(c_0) = a_0$;
$a(c) = f_{c}^{p}(a(c))$; and $|(f_{c}^p)'(a(c))| >1$.
We let $\lam(c):= (f_{c}^p)'(a(c))$ and 
$\lam_0:= \lam(c_0)$. 

Our key lemma is the following.

\begin{lem}\label{lem_key}
Suppose that $c_0 \in \M$ is a Misiurewicz parameter as above. 
For $k \in \N$, set $\rho_k:= 1/(f_{c_0}^{l + kp})'(c_0)$.
Then we have the following.
\begin{enumerate}[\rm (1)]
\item
The function $\phi_k(w) = f_{c_0}^{l + kp}(c_0 + \rho_k w)$ 
converges to a non-constant entire function 
$\phi:\C \to \C$ as $k \to \infty$ uniformly on any compact sets.
\item
There exists a constant $Q \neq 0$ such that the function 
$$
\Phi_k(w) \dee 
f_{c_0 + Q\rho_k w}^{l + kp}(c_0 + Q\rho_k w)
$$
converges to the same function $\phi(w)$ as $k \to \infty$ 
uniformly on compact sets of $\C$.
\end{enumerate}
\end{lem}

\begin{pf}
It is well-known 
that the sequence of (polynomial) functions 
$$
w ~\longmapsto~ f_{c_0}^{kp}\paren{a_0 + \frac{w}{\lam_0^k}} 
~~~~~~~(k \in \N)
$$
converges to a non-constant entire function $\phi(w)$ uniformly on compact sets of $\C$.
(See \thmref{thm_Poincare_func} in Appendix.
Such a $\phi$ is called a {\it Poincar\'e function}.
Indeed, $\phi$ satisfies the functional equation
$
\phi(\lam_0w ) = f_{c_0}^{p} \cc \phi(w),
$
but we will not use it.)
%Such a $\phi$ is called a {\it Poincar\' e function}.)
Note that this function satisfies
$\phi(0) = a_0$ and $\phi'(0) = 1$.

Now let us show (1): 
set $A_0 := (f^l_{c_0})'(c_0)$, 
where $A_0\neq 0$ since otherwise $c_0$ is strictly periodic.
We also have $(f^{l+kp}_{c_0})'(c_0)= A_0 \lam_0^k=1/\rho_k$.
For sufficiently small $t \in \C$, we have the expansion
$$ 
f^l_{c_0}(c_0 + t) 
\ee a_0 + A_0 \cdot t + o(t).
$$ 
Fix an arbitrarily large compact set $E \subset \C$ 
and take any $w \in E$.
Then by setting $t= w/(A_0\lam_0^k)$,
$$ 
f_{c_0}^{kp}\paren{a_0 + \frac{w}{\lam_0^k}} 
~\sim~ 
f^{l+kp}_{c_0}
\paren{c_0 + \frac{w}{A_0\lam_0^k} +  o(\lam_0^{-k})}
~\sim~
f_{c_0}^{l+kp}\paren{c_0 + \rho_k w} 
~ = ~ \phi_k(w)
$$ 
when $k \to \infty$. 
(Here by $A_k(w) \sim B_k(w)$ we mean $A_k(w)-B_k(w) \to 0$ uniformly on $E$ as $k \to \infty$.)
Hence we have
$ 
\phi(w) =  \lim_{k \to 0} \phi_k(w)
$
on any compact sets.
Note that $\phi$ has no poles, since each $\phi_k$ is entire.

Next we show (2): 
suppose that $Q \in \Cstar$ is a constant 
and set $c = c(w):= c_0 + Q \rho_k w$.
We also set $\Phi_k(w) := f_{c}^{l +  kp}(c)$ 
and $b(c):= f_c^{l}(c)$.
Recall that $a(c)$ denotes a repelling periodic point 
(of $f_c$) of period $p$ with $a(c_0) = a_0$, 
and $\lam(c)$ denotes its multiplier.

Then \thmref{thm_Poincare_func} in Appendix again implies that 
the sequence of functions 
$\phi_{k}^{c}(w):= f_c^{kp}(a(c) + w/\lam(c)^{k})$ 
converges to an entire function 
$\phi^c(w)$ uniformly on compact sets.
In particular, by the proof of 
\thmref{thm_Poincare_func}, 
it is not difficult to check that 
the function $c \mapsto \phi^c(w)$ is holomorphic 
near $c = c_0$ when we fix a $w \in \C$.

As in Tan's original proof, we employ a theorem on transversality 
by Douady and Hubbard \cite[Lemma 1, p.333]{DH2}:
There exists a $B_0 \neq 0$ such that 
$$
b(c)-a(c) = B_0(c-c_0) +o(c-c_0).
$$
Hence for $c = c_0 + Q\rho_k w$ (taking $w$ in a compact set), we have 
$$
b(c) \ee a(c) + B_0Q \rho_k w +o(\rho_k) 
\ee a(c) + 
\frac{B_0Q}{A_0} 
 \cdot 
\frac{ \lam(c)^{k} }{ \lam_0^{k}} 
 \cdot \frac{w}{\lam(c)^{k}}
+o(\rho_k). 
$$
Set $Q:= A_0/B_0$. 
Since $\lam(c)$ is a holomorphic function of $c$
and thus $\lam(c) = \lam_0 + O(c-c_0)$,
we have $|\lam(c)/\lam_0-1| = O(c-c_0)$.
This implies that 
$$
\log \frac{ \lam(c)^{k} }{ \lam_0^{k}} 
 = k \cdot O(c-c_0) 
 = O\paren{\frac{k}{\lam_0^{k}}}
  \to 0~~(k \to \infty)
$$
for $c = c_0 + Q\rho_k w = c_0 + O(\lam_0^{-k})$.
Since 
$$
\Phi_k(w) =f_c^{kp}(b(c))\ee 
f_c^{kp}\paren{a(c) +  \frac{w}{\lam(c)^{k}}+o(\rho_k)}
$$ 
and 
$\phi^c(w) \to \phi(w)$ as $c \to c_0$
(uniformly for $w$ in a compact set), 
we conclude that 
$$
\lim_{k \to \infty} \Phi_k(w) 
\ee \phi(w),
$$ 
where the convergence is uniform on any compact sets.
\QED
\end{pf}

\parag{Remark.}
Lemma 1 implies that $c_0 \in J_{c_0}=\partial K_{c_0}$ and $c_0 \in \partial \M$. Indeed, we can find a $w \in \C$ such that $|\phi(w)|>2$ 
and hence $|\phi_k(w)|>2$ for sufficiently large $k$. 
Equivalently, we have $c_0+\rho_k w \notin K_{c_0}$ 
for sufficiently large $k$, where $c_0+\rho_k w$ 
tends to $c_0$ as $k \to \infty$.  
Since $c_0 \in K_{c_0}$ by definition, we have $c_0 \in J_{c_0}$.
The proof for $c_0 \in \partial \M$ is analogous.

\parag{The Hausdorff topology.}
Let us briefly recall the \text{\it Hausdorff topology} of the set of non-empty compact sets $\mathrm{Comp}^\ast(\C)$ of $\C$. 
For a sequence $\{K_k\}_{k\in \N} \subset \mathrm{Comp}^\ast(\C)$, we say $K_k$ converges to $K \in \mathrm{Comp}^\ast(\C)$ as $k \to \infty$ if for any $\e>0$, 
there exists $k_0 \in \N$ such that 
$K \subset \mathrm{N}_\e(K_k)$ and $K_k \subset \mathrm{N}_\e(K)$
for any $k \ge k_0$, 
where $\mathrm{N}_\e(\cdot)$ is the open $\e$ neighborhood in $\C$.

Let $\D(r): = \braces{z \in \C \st |z| < r}$.
For a compact set $K$ in $\C$, 
let $[K]_r$ denote the set $(K \cap \D(r))\cup \partial \D(r)$.
For $a \in \Cstar$ and $b \in \C$, 
let $a(K-b)$ denote the set of $a(z-b)$ with $z \in K.$

\parag{Similarity.}
Let $c_0$ be a Misiurewicz parameter.
Now we state our version of Tan's similarity theorem.

\begin{thm}[Similarity between $\M$ and $J$]\label{thm_MJ}
There exist a non-constant entire function $\phi$ on $\C$, 
a sequence $\rho_k \to 0$,
and a constant $q \neq 0$
such that if we set $\cJ: = \phi^{-1}(J_{c_0}) \subset \C$,
then for any large constant $r>0$, we have
\begin{enumerate}[\rm (a)]
\item 
$\gauss{\rho_k^{-1}(J_{c_0}-c_0)}_r  \to [\cJ]_r$, and 
\item 
$\gauss{\rho_k^{-1}q(\M-c_0)}_r  \to [\cJ]_r$%~~(k \to \infty)$
\end{enumerate}
as $k \to \infty$ in the Hausdorff topology. 
\end{thm}
%Actually $\phi$, $\rho_k~(k \in \N)$ 
%and $1/q = Q$ are all given in (the proof of) \lemref{lem_key}.

\paragraph{Proof of (a).}
Let $\phi_k,~\phi$, and $\rho_k = 1/(A_0\lam_0^k)$ be as given in the 
proof of \lemref{lem_key}.
Since $f_{c_0}^{n}(J_{c_0}) = J_{c_0}$, we have $\gauss{\rho_k^{-1}(J_{c_0}-c_0)}_r = \gauss{\phi_k^{-1}(J_{c_0})}_r$. 
By $[\cJ]_r = \gauss{\phi^{-1}(J_{c_0})}_r$ and uniform convergence of $\phi_k \to \phi$ on $\overline{\D(r)}$, the claim easily follows.

\paragraph{Proof of (b).}
Let $q := 1/Q$ (where $Q \neq 0$ is defined in the proof of \lemref{lem_key}) and $\cM_k:= \rho_k^{-1}q(\M-c_0)$. 
Fix any $\e>0$. Since the set $\overline{\D(r)}-\mathrm{N}_\e(\cJ)$ is compact, 
there exists an $N = N(\e)$ such that $|f_{c_0}^N \cc \phi(w)| >2$ for any $w \in \overline{\D(r)}-\mathrm{N}_\e(\cJ)$. 
By uniform convergence of 
$\Phi_k(w) = f^{l + kp}_{c_0 + Q\rho_k w}(c_0 + Q\rho_k w)$
to $\phi(w)$ on compact sets in $\C$ (\lemref{lem_key}), 
we have 
$$
|f_{c_0 + Q \rho_k w}^{(l + kp) + N}(c_0 + Q \rho_k w)| ~>~2
$$
for all sufficiently large $k$. This implies that $c_0 + Q\rho_k w \notin \M$, equivalently, $w \notin \cM_k$. 
Hence we have 
$$
\gauss{\cM_k}_r ~\subset~ \mathrm{N}_\e([\cJ]_r).
$$ 

Next we show the opposite inclusion $[\cJ]_r\subset \mathrm{N}_\e(\gauss{\cM_k}_r)$ 
for $k$ large enough.
Let us approximate $[\cJ]_r$ by a finite subset $E$ of $[\cJ]_r$ such that the $\e/2$ neighborhood of $E$ covers $[\cJ]_r$. 
Now it is enough to prove that for any $w_0 \in E$, there exists a sequence $w_k \in [\cM_k]_r$ such that $|w_0-w_k| < \e/2$ for sufficiently large $k$.

Let $\Delta$ be a disk of radius $\e/2$ centered at $w_0$. 
When $\Delta \cap \del \D(r) \neq \emptyset$, we can take such a $w_k$ in $\del \D(r)$. Hence we may assume that $\Delta \subset \D(r)$.

Since $\phi(w_0) \in J_{c_0}$ and repelling cycles are dense in $J_{c_0}$
(see \cite{Sch} and \cite{Mi2}. See also the remark below), we can choose a $w_0'$ such that $\phi(w_0')$ is a repelling periodic point of some period $m$ and $|w_0-w_0'|< \e/4$. 
This implies that the function $\chi:w \mapsto f_{c_0}^m(\phi(w))-\phi(w)$ has a zero at $w = w_0'$. 

Let us consider the function 
$\chi_k: w \mapsto f_{c_0 + Q \rho_k w}^m(\Phi_k(w))-\Phi_k(w)$, 
where $\Phi_k(w) = f_{c_0 + Q \rho_k w}^{l + kp}(c_0 + Q \rho_k w)$
as in \lemref{lem_key}. 
By the Hurwitz theorem and uniform convergence of $\Phi_k$ to $\phi$ on compact sets of $\C$, $\chi_k$ has a zero $w_k$ in $\Delta$ and $|w_k - w_0'| < \e /4$ for all sufficiently large $k$. In particular, $c_k:= c_0 + Q \rho_k w_k$ satisfies $f_{c_k}^{m + (l + kp)}(c_k) = f_{c_k}^{l + kp}(c_k)$ and thus $c_k \in \M$. Hence we have a desired $w_k \in \cM_k$ with $|w_k -w_0|<\e/2$. 
\hfill \QED

\medbreak

\parag{Example (Calculation of $Q$).}
When $c_0=-2$ (hence $l=p=1$), 
$A_0=f_{c_0}'(-2)=-4$
and $\lam_0 = f_{c_0}'(f_{c_0}(-2))=4.$
Since $a(c)^2+c=a(c)$ we find $da(c)/dc=-1/(2 a(c)-1)$.
Moreover, $db(c)/dc=(d/dc)(c^2+c)=2c+1$. Hence for $c=c_0=-2$,
we have $B_0=-3-(-1/3)=-8/3$ and the constant $Q$ is 
$A_0/B_0=3/2$.

\parag{Remarks.}
\begin{itemize}
\item
Note that in the proof of \thmref{thm_MJ}, $\phi(w_0')$ need not be repelling.
We only need the density of periodic points in the Julia set,
which is an easy consequence of Montel's theorem.
See \cite[p.157]{Mi2}. 

\item
A similar proof can be applied to semi-hyperbolic parameters
(i.e., critically non-recurrent and no parabolic cycle) in $\partial \M$
 \cite[Th\'eor\`eme 2.2]{Ka}, and 
 to the unicritical family $\braces{z \mapsto z^d+c \st c\in \C}$
 with $d \ge 2$.  
This gives an alternative proof of Rivera-Letelier's 
extension of Tan's theorem \cite{RL}. 
\end{itemize}

\section{Zalcman's lemma and Zalcman functions}
The key to the proof above is \lemref{lem_key}
that bridges the dynamical and parameter planes 
by one entire function $\phi$.
Let us characterize this property in a generalized setting
by means of {\it Zalcman's lemma},
which gives a precise condition for non-normality:

\parag{Zalcman's lemma {\cite{Za, Za2}}.}
{\it 
Let $D$ be a domain in the complex plane $\C$ and 
$\cF$ a family of meromorphic functions on $D$.
The family $\cF$ is not normal on any neighborhood of $z_0 \in D$
if and only if
there exist sequences
$F_k \in  \cF$,
%$n_k \in \N$ with $n_k \to \infty$; 
$\rho_k \in \Cstar$ with $\rho_k \to 0$; and
$z_k \in D$ with $z_k \to z_0$ such that 
the function $\phi_k(w) = F_{k}(z_k + \rho_k w)$ 
converges to a non-constant meromorphic function
$\phi:\C \to \Chat = \C \cup \braces{\infty}$ 
uniformly on compact subsets in $\C$.
}

\parag{A universal setting.}
Let $\cUbar$ be the space of meromorphic functions on $\C$
 with a topology induced by uniform convergence on compact subsets
in the spherical metric. 
Let $\cU \subset \cUbar$ be the set 
of non-constant meromorphic functions on $\C$
and $\Aff \subset \cU$ 
the set of complex affine maps.
One can easily check that Zalcman's lemma above can be translated as follows:
{\it 
A family of meromorphic functions $\cF$ on the domain $D$
is not normal in any neighborhood of $z_0 \in D$
if and only if there exist $A_k \in \Aff$ and $f_k \in \cF$ $(k \in \N)$
such that as $k \to \infty$, 
\begin{enumerate}
\item
$A_k$ converges to a constant function $z_0$ in $\cUbar$, and
\item
$f_k \cc A_k$ converges to some $\phi$ in $\cU$.
\end{enumerate}
}
We denote the set of all possible limit function $\phi \in \cU$ 
of this form by
$
\cZ(\cF,z_0),
$
and we say $\phi \in \cZ(\cF,z_0)$ is a {\it Zalcman function
of the family $\cF$ at $z_0$}. 
If $\cF$ is normal on a neighborhood of $z_0$, we formally set
$
\cZ(\cF,z_0):=\emptyset.
$
Now the union
$$
\cZ(\cF):=\bigcup_{z_0 \in D} \cZ(\cF,z_0) \subset \cU
$$
is the set of Zalcman functions of the family $\cF$.

\if0
The original statement to this lemma is straightforward.
To check the converse, let $A_k(w):= z_k + \rho_k w \in \Aff$
and $f_k \cc A_k$ converges to $\phi \in \cU$.
Then $z_k \to z_0$ and $\rho_k \to 0$ as $k \to \infty$.
By extracting a subsequence, we may assume that 
$\arg \rho_k$ converges to some $\theta$
and thus $f_k(z_k+|\rho_k|w)$ converges to $\phi(e^{ -i\theta} w)$
in $\cU$.
\fi

\parag{Dynamical and parametric Zalcman functions.}
We want to give a new perspective on the asymptotic similarity between 
the Mandelbrot set $\M$ and the Julia set $J_c$ in terms of Zalcman's lemma.
Let us recall the following facts.

\begin{prop}[$J$ and $\partial \M$ as non-normality loci]
\label{prop_non-normal_MJ}
For the quadratic family $f_c(z)=z^2+c~(c \in \C)$, 
the Julia sets and the boundary of the Mandelbrot set
are characterized as follows:
\begin{enumerate}[\rm (1)]
\item
For each $c \in \C$, the Julia set $J_c=\partial K_c$ 
is the set of points where the family 
$\braces{z \mapsto f_c^n(z)}_{n \ge 0}$ is not normal.
\item
The boundary $\partial \M$ of the Mandelbrot set is the set 
of points where the family 
$\braces{c \mapsto f_c^n(c)}_{n \ge 0}$ is not normal.
\end{enumerate}
\end{prop}
The proof is done by an elementary equicontinuity argument. 
%For the Julia set case, see \cite{Mi}?. 
For the Mandelbrot set case, see \cite[Theorem 4.6]{Mc}.

\parag{Dynamical Zalcman functions.}
Let $\cF_c$ denote the family $\braces{z \mapsto f_c^n(z)}_{n \ge 0}$ 
of polynomial functions on $\C$
for each parameter $c \in \C$.
By \propref{prop_non-normal_MJ}, we can apply Zalcman's lemma
to this family and obtain the sets 
$$
\cZ_c:=\cZ(\cF_c)
\quad\text{and}\quad 
\cZ_c(z_0):=\cZ(\cF_c, z_0) 
$$
of {\it dynamical Zalcman functions} of $f_c$ (for each $z_0 \in J_c$).
These sets have a good invariance with respect to 
the operations `$f_c \cc$' and `$\cc \Aff$',
 as shown by Steinmetz \cite[Theorems 1 and 2]{St}:

\begin{prop}[Invariance]\label{prop_St}
For each $z_0 \in J_c$, the family $\cZ_c(z_0)$ satisfies 
$$
f_c\cc \cZ_c(z_0) \ee \cZ_c(z_0) \ee \cZ_c(z_0)\cc \Aff.
$$
More precisely, 
\begin{enumerate}[\rm (1)]
\item
If $\phi \in \cZ_c(z_0)$ then $f_c\cc \phi \in \cZ_c(z_0)$
and $\phi = f_c \cc \phi_1$ for some $\phi_1 \in \cZ_c(z_0)$. 
\item 
For any $A \in \Aff$ and $\phi \in \cZ_c(z_0)$, 
we have $\phi \cc A^{\pm 1}  \in \cZ_c(z_0)$. 
\end{enumerate}
\end{prop}
Note that the universal space $\cU$ only satisfies 
$f_c \cc \cU \subset \cU$ and $\cU \cc \Aff = \cU$.

\propref{prop_St} implies that we also have 
$
f_c \cc \cZ_c = \cZ_c =\cZ_c \cc \Aff.
$ 
However, the equality $\cZ_c=\cZ_c(z_0)$ holds for any $z_0 \in J_{c}$ 
in most cases. To see this, we introduce some terminology: 
The {\it univalent grand orbit} $UGO(z_0)$ of $z_0 \in \C$
is the set of $\zeta$ such that $f_c^m(z_0)=f_c^n(\zeta)$
for some $m,\,n \in \N$ and 
there is a univalent branch $g$ of $f_c^{-n}\cc f_c^m$ 
in a neighborhood of $z_0$ with $g(z_0)=\zeta$.
The {\it postcritical set} $P_c$ of $f_c$ is the closure of the orbit 
$\braces{c, f_c(c), f_c^2(c), \cdots}$.
(Note that $c$ is a unique critical value of $f_c$.)
We say {\it $f_c$ satisfies $(\ast)$-condition} if
\begin{quote}
$(\ast)$ : For any $z_0 \in J_c$, there exists a $\zeta \in UGO(z_0)-P_c$.
\end{quote}
Now we claim the following:
\begin{thm}\label{thm_coincidence}
If $f_c$ satisfies $(\ast)$-condition, 
then $\cZ_c=\cZ_c(z_0)$ for any $z_0 \in J_{c}$.
Moreover, the set of such $c$ contains $\C-\partial \M$, which is an open and dense subset of $\C$, and the Misiurewicz parameters in $\partial \M$ except $c_0=-2$. 
\end{thm}

\begin{pf}
The first claim of the theorem is an immediate corollary 
of \cite[Th\'eor\`eme 3.8]{Ka}. 
For the second claim, note that $\partial \M$ consists 
of the parameters $c$ 
for which $f_c$ is expanding (hyperbolic) or infinitely renormalizable
(see \cite[\S 4 and \S 8]{Mc}).
Since expanding $f_c$ does not contain any critical point in the Julia set $J_c$, it satisfies $(\ast)$-condition. 
Any infinitely renormalizable $f_c$ satisfies $(\ast)$-condition by \cite[Proposition 3.7]{Ka}.

Now suppose that $c$ is a Misiurewicz parameter,
satisfying $f_{c}^l(c)=f_{c}^{l+p}(c)$ with minimal $l,\,p \ge 1$.
Then $z_0 \in P_c$ if and only if 
$z_0=f_c^k(c)$ for some $0 \le k \le l+p-1$.
For such a $z_0 \in P_c$, 
when $l \ge 2$ or $p \ge 2$, 
one can easily find some $\zeta \in UGO(z_0)-P_c$.
Otherwise $l=p=1$, equivalent to $c=-2$.
In this case we have $P_{-2}=\{-2,\,2\}=UGO(-2)=UGO(2)$.

If $z_0 \notin P_c$, $z_0$ itself is an element of $UGO(z_0)-P_c$.
Hence $f_c$ for the Misiurewicz parameter $c$ satisfies $(\ast)$-condition unless $c=-2$.
\QED\end{pf}

\parag{Remark.}
When $f_c$ satisfies $(\ast)$-condition, 
the dynamical Zalcman function $\cZ_c$ can be an alternative ingredient 
for the Lyubich-Minsky Riemann surface lamination for $f_c$. 
See \cite[\S 3]{Ka} for more details, 
and see \cite{LM} for Lyubich and Minsky's 
lamination theory in complex dynamics.

\parag{Parametric Zalcman functions.}
Let $\cQ$ denote the family $\braces{c \mapsto f_c^n(c)}_{n \ge 0}$ 
of polynomial functions on $\C$.
Again by \propref{prop_non-normal_MJ}, 
we can apply Zalcman's lemma
to this family and obtain the sets 
$$
\cP:=\cZ(\cQ)
\quad\text{and}\quad 
\cP(c_0):=\cZ(\cQ, c_0) 
$$
of {\it parametric Zalcman functions} of the quadratic 
family $\braces{f_c(z)=z^2+c}_{c \in \C}$
 (for each $c_0 \in \partial \M$).
These sets have weaker invariance,
which is also pointed out by Steinmetz \cite[Remark in \S 1]{St}:

\begin{prop}[Invariance for $\cP$]\label{prop_St_para}
For each $c_0 \in \partial \M$, the family $\cP(c_0)$ satisfies 
$$
f_{c_0}\cc \cP(c_0) \ee \cP(c_0) \ee \cP(c_0)\cc \Aff.
$$
\end{prop}
Hence we only have $\cP = \cP \cc \Aff$ for the total space
$\cP$ of the parametric Zalcman functions.

In a forthcoming paper we will present a general account on 
dynamical and parametric 
Zalcman functions for families of rational functions
parametrized by Riemann surfaces. 

\parag{Dynamical-parametric intersection and similarity.}
By \propref{prop_St} and \propref{prop_St_para} above,
if $c_0 \in J_{c_0}$ and $c_0 \in \partial \M$,
then $\cZ_{c_0}(c_0)$ and $\cP(c_0)$ 
exhibit the same invariance in the universal space $\cU$. 
Hence one might expect that there exists some 
$\phi \in \cZ_{c_0}(c_0) \cap \cP(c_0)$
when $c_0 \in J_{c_0} \cap \partial \M$.
Indeed, the existence of such an intersection implies 
asymptotic similarity between $J_{c_0}$ and $\M$ at $c_0$:

\begin{thm}[Intersection implies similarity]\label{thm_intersection_implies_similarity}
Suppose that $\cZ_{c_0}(c_0) \cap \cP(c_0) \neq \emptyset$
for some $c_0 \in \partial \M$. 
More precisely, there exist sequences of affine maps 
$A_k,\, B_k \in \Aff$ 
and positive integers $m_k, \, n_k  \in \N$ 
such that as $k \to \infty$ we have 
\begin{enumerate}[\rm (1)]
\item
Both $A_k$ and $B_k$ converge to the same constant map $c_0$ in $\cUbar$;
and
\item
Both 
$f_{c_0}^{n_k} (A_k(w))$ and 
$f_{B_k(w)}^{n_k} ( B_k(w))$ 
converge to the same entire function $\phi(w)$ 
in $\cU$.
\end{enumerate}
Suppose in addition that $c_0 \in J_{c_0}$ 
and let $\cJ: = \phi^{-1}(J_{c_0}) \subset \C$. 
Then for any large constant $r>0$, we have
\begin{enumerate}[\rm (a)]
\item 
$\gauss{A_k^{-1}(J_{c_0})}_r  \to [\cJ]_r$, and 
\item 
$\gauss{B_k^{-1}(\M)}_r  \to [\cJ]_r$%~~(k \to \infty)$
\end{enumerate}
as $k \to \infty$ in the Hausdorff topology. 
\end{thm}

\parag{Proof.}
Condition (1) implies that if we write
$A_k(w)=c_k + \rho_k w$ and $B_k(w)=c_k' + \rho_k' w$
then $c_k, c_k' \to c_0$ and $\rho_k, \rho_k' \to 0$
as $k \to \infty$. Then the proof is only a slight modification 
of that of \thmref{thm_MJ}.
\QED

\parag{Intersection at Misiurewicz parameters.}
Suppose that $c_0$ is a Misiurewicz parameter.
(Hence $c_0 \in J_{c_0}$ and $c_0 \in \partial \M$ 
by the remark below the proof of Lemma 1.)
By the construction of the entire function $\phi$ in Lemma 1,
we obviously have $\phi \in \cZ_{c_0}(c_0)$ and 
$\phi \in \cP(c_0)$. 
Thus we conclude the following:

\begin{thm}[Intersection]\label{thm_intersection}
For any Misiurewicz parameter $c_0$, 
$\cZ_{c_0}(c_0)$ and $\cP(c_0)$ share at least one element $\phi \in \cU$.
\end{thm}

Note that the set of Misiurewicz parameters is a dense subset of
$\partial \M$. 
(This can be shown by a standard normal family argument. 
See Levin \cite{Le} or \cite[Th\'eor\`eme 1.1 (1)]{Ka}.)

In \cite{Ka} the author proved Lemma 1 for a wider class of parameters 
 in $\partial \M$, called {\it semi-hyperbolic parameters}. 
Shishikura proved in \cite{Sh} that the set of semi-hyperbolic parameters is a dense subset of $\partial \M$ of Hausdorff dimension 2. 
Hence we also have a generalization of \thmref{thm_intersection} for 
semi-hyperbolic parameters.

\parag{Question.} 
Besides the example of Lemma 1 and its generalization to semi-hyperbolic parameters, are there any other intersections of the sets of dynamical and parametric 
Zalcman functions?

\section{Similarity between $\T$ and $J$}
We apply the arguments in Section \ref{sec_MJ} 
to the antiholomorphic quadratic family
and we will show an analogous result to Tan's theorem.
For the theory of antiholomorphic quadratic family
and the tricorn (as a counterpart of the Mandelbrot set), 
see 
\cite{CHRS},
\cite{Mi1},
\cite{N},
\cite{NS},
\cite{MNS},
\cite{HS},
\cite{I}, and 
\cite{IM}
for example.

\parag{The tricorn and the Julia sets.}
Let us consider the antiholomorphic quadratic family 
$$
\braces{g_c(z) = \zbar^{\,2} + c \st c \in \C}.
$$ 
The {\it tricorn} $\T$ is the set of $c \in \C$ 
such that the sequence $\braces{g_c^n(c)}_{n \in \N}$ is bounded. 
For each $c \in \C$, the \text{\it filled Julia set} $K_c^\ast$ 
is the set of $z \in \C$ such that the sequence 
$\braces{g_c^n(z)}_{n \in \N}$ is bounded.
One can easily check that
\begin{itemize}
\item
$c \notin \T$ if and only if 
$|g_c^n(c)| > 2$ for some $n \in \N$; and
\item
for each $c \in \T$, $z \notin K_c^\ast$ if and only if 
$|g_c^n(z)| > 2$ for some $n \in \N$.
\end{itemize}
The {\it Julia set} $J_c^\ast$ is the boundary of $K_c^\ast$.
Note that all $\T, K_c^\ast$, and $J_c^\ast$ 
are non-empty compact sets.

An intriguing property of $\T$ is that 
$\T=\omega \T=\omega^2 \T$ for $\omega=e^{2 \pi i /3}$. 
One can also check that $\T \cap \R= \M \cap \R=[-2,1/4]$, 
and we have $K_c=K_c^\ast$ and $J_c=J_c^\ast$ for real $c$.

%%%%%%----------table of pictures
\fboxsep=0pt
\fboxrule=.5pt
\begin{figure}[htbp]
\begin{center} 
\phantom{(T)}
\fbox{\includegraphics[width=.45\textwidth, bb = 0 0 800 800]{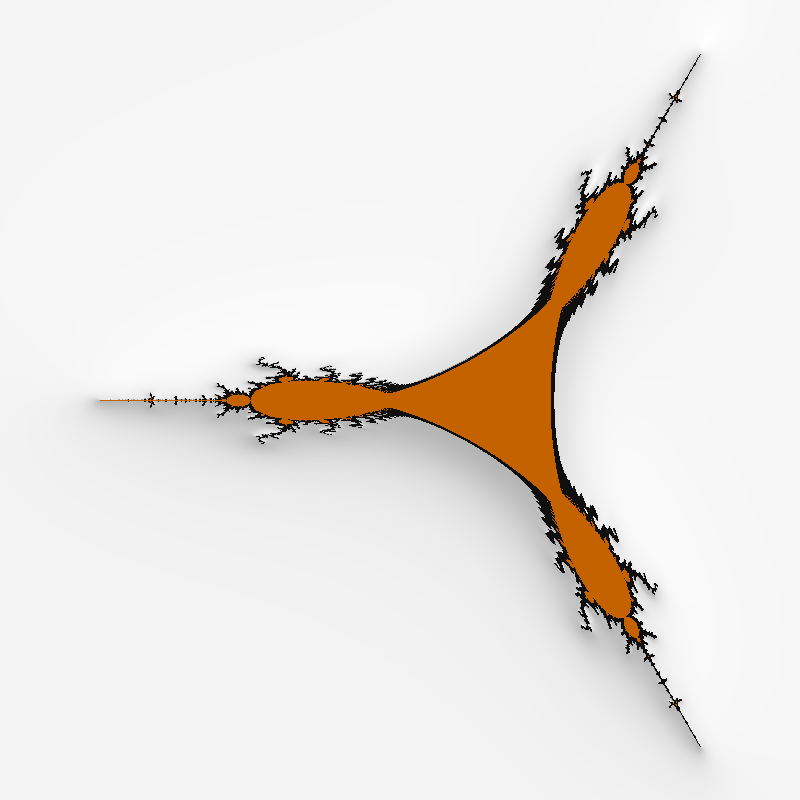}}
(T)\\[.5em]
\fbox{\includegraphics[width=.18\textwidth, bb = 0 0 600 600]{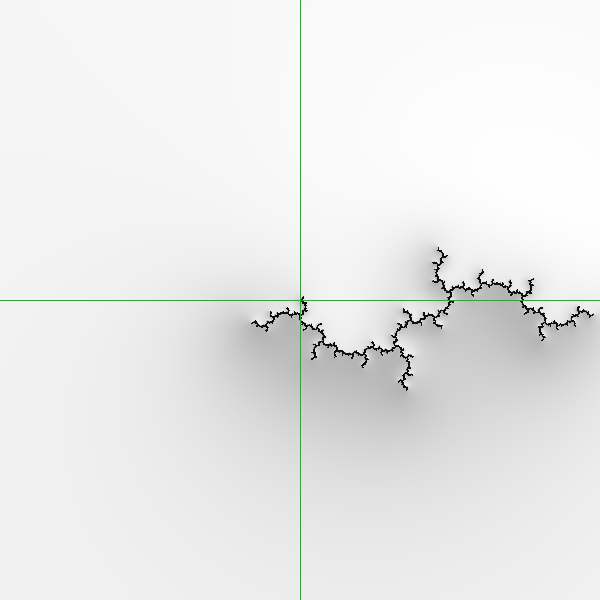}}
\fbox{\includegraphics[width=.18\textwidth, bb = 0 0 600 600]{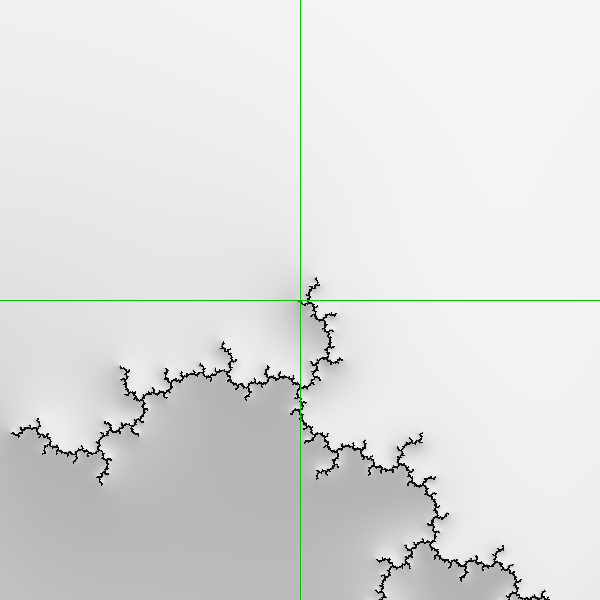}}
\fbox{\includegraphics[width=.18\textwidth, bb = 0 0 600 600]{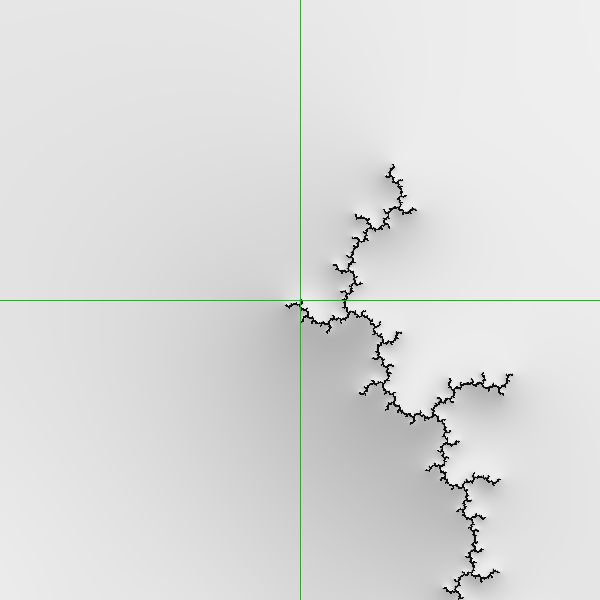}}
\fbox{\includegraphics[width=.18\textwidth, bb = 0 0 600 600]{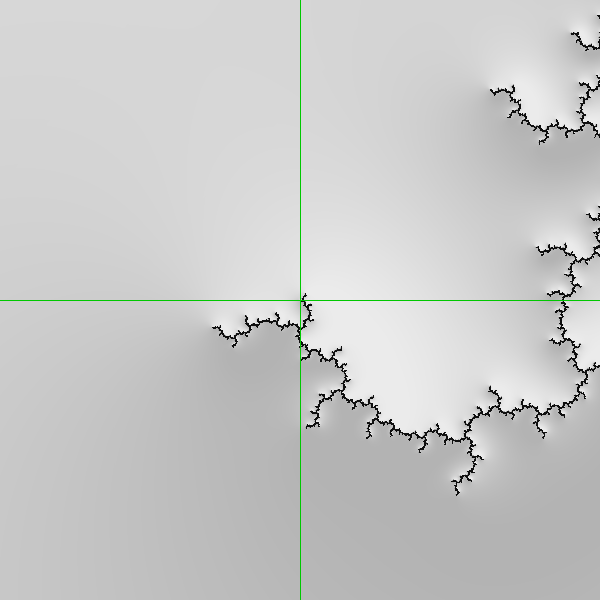}}
\fbox{\includegraphics[width=.18\textwidth, bb = 0 0 600 600]{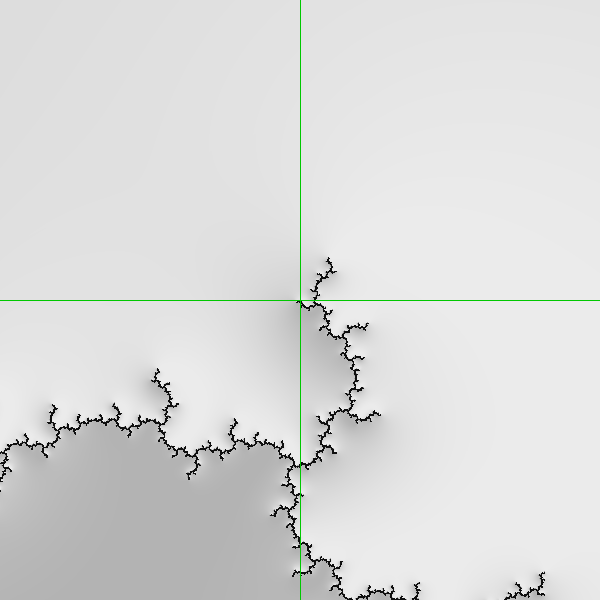}}
\\[.5em]
\fbox{\includegraphics[width=.18\textwidth, bb = 0 0 600 600]{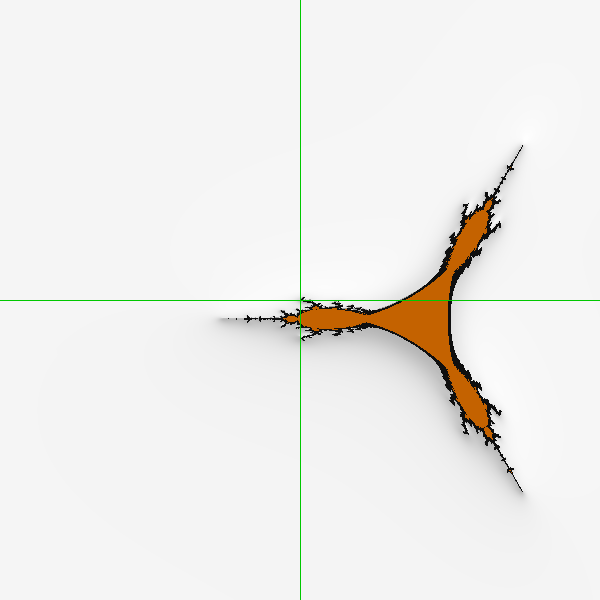}}
\fbox{\includegraphics[width=.18\textwidth, bb = 0 0 600 600]{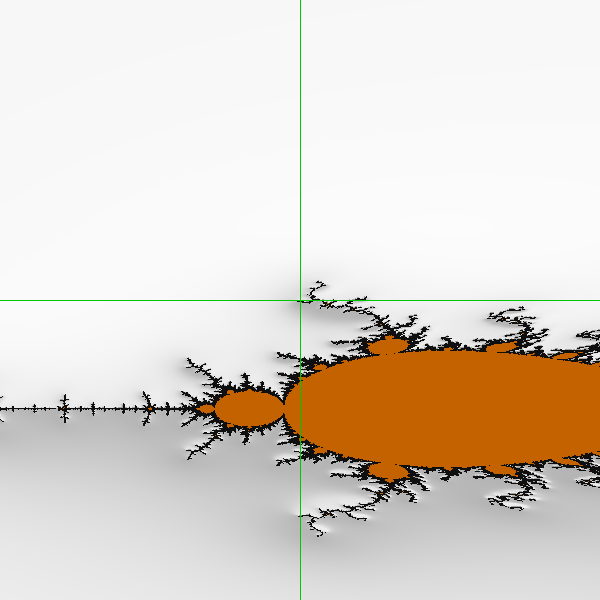}}
\fbox{\includegraphics[width=.18\textwidth, bb = 0 0 600 600]{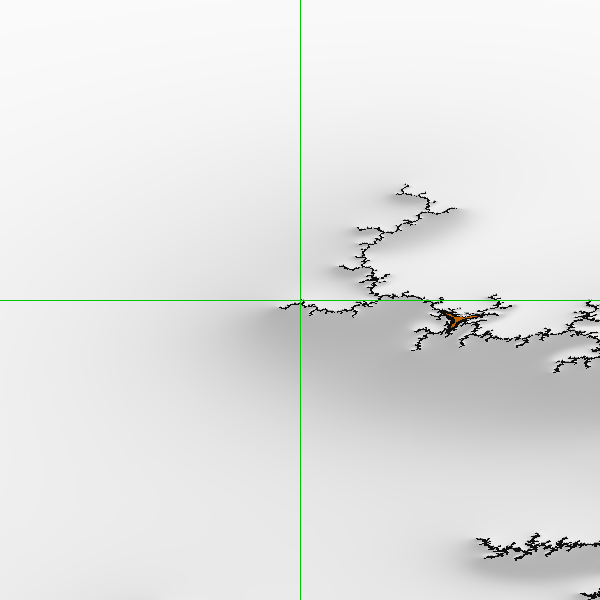}}
\fbox{\includegraphics[width=.18\textwidth, bb = 0 0 600 600]{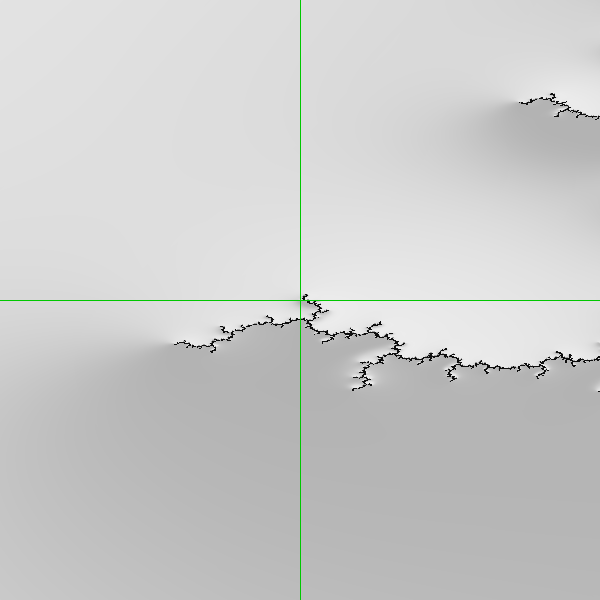}}
\fbox{\includegraphics[width=.18\textwidth, bb = 0 0 600 600]{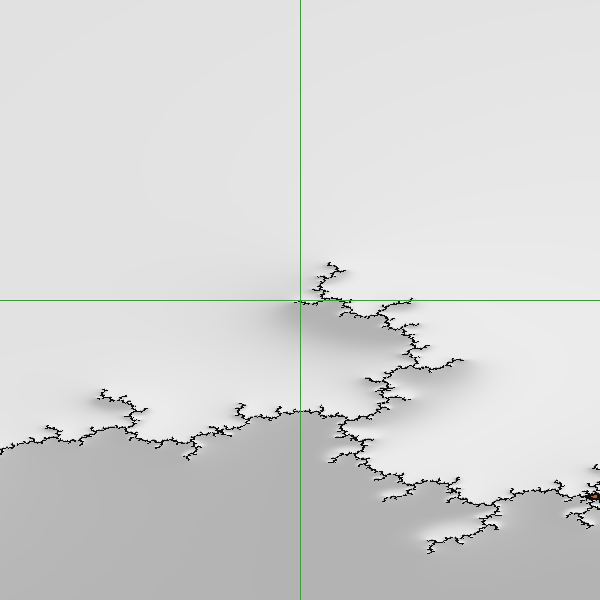}}
\\[.5em]
(JT1)\\[.5em]%:Center: $0.0+1.0 i$, square width: from 5.0 to 0.01. \\[.5em]
\fbox{\includegraphics[width=.18\textwidth, bb = 0 0 600 600]{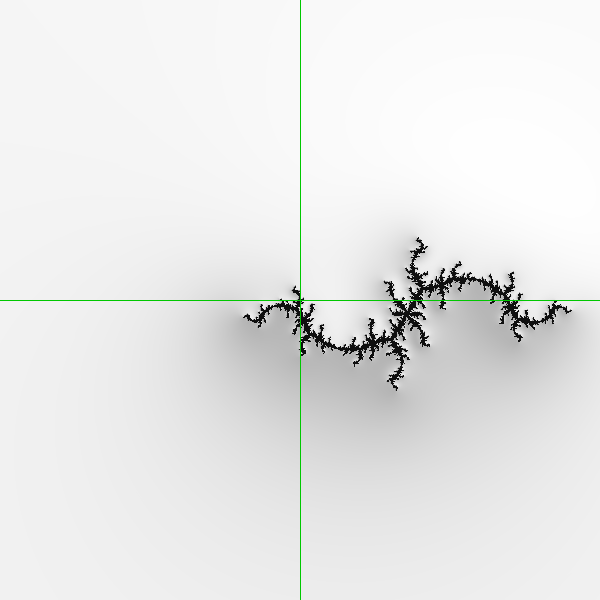}}
\fbox{\includegraphics[width=.18\textwidth, bb = 0 0 600 600]{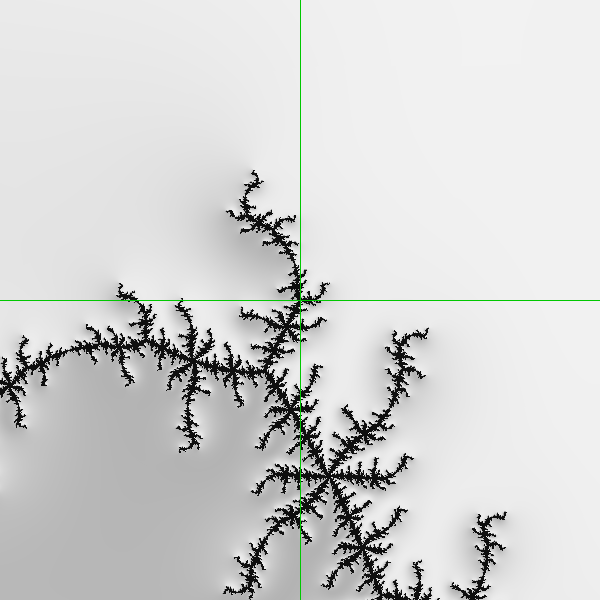}}
\fbox{\includegraphics[width=.18\textwidth, bb = 0 0 600 600]{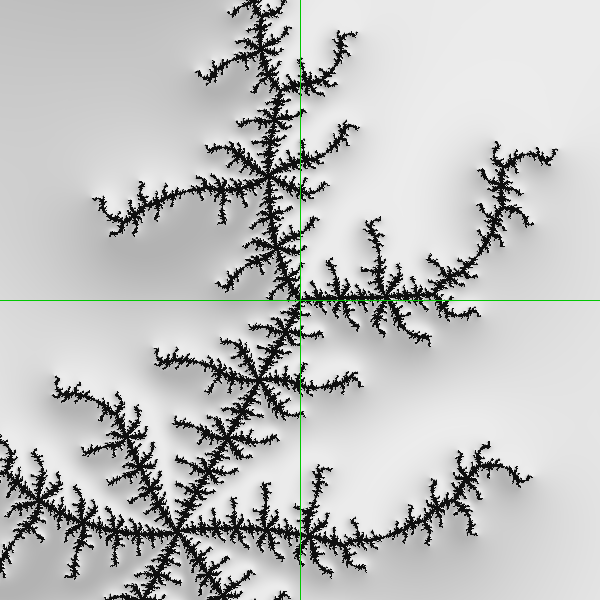}}
\fbox{\includegraphics[width=.18\textwidth, bb = 0 0 600 600]{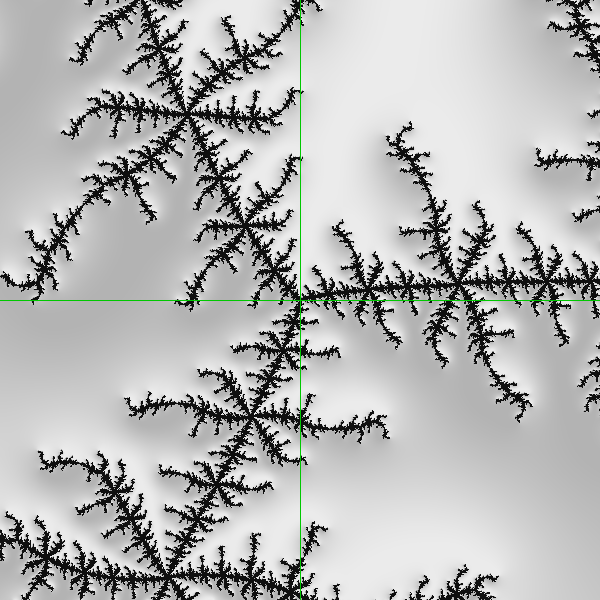}}
\fbox{\includegraphics[width=.18\textwidth, bb = 0 0 600 600]{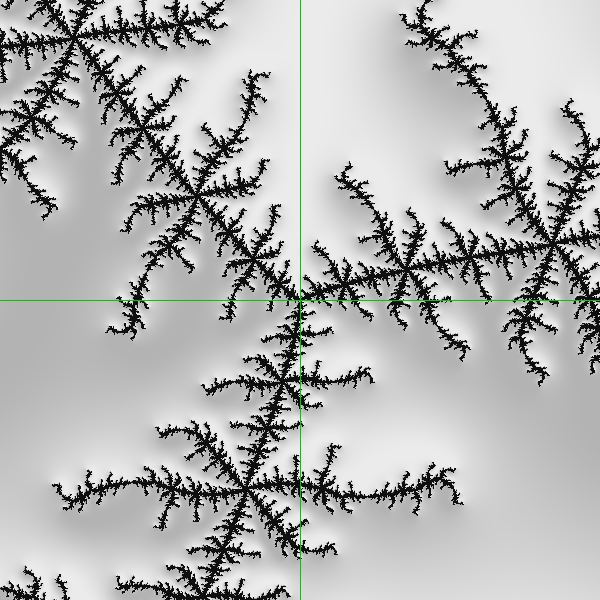}}
\\[.5em]
\fbox{\includegraphics[width=.18\textwidth, bb = 0 0 600 600]{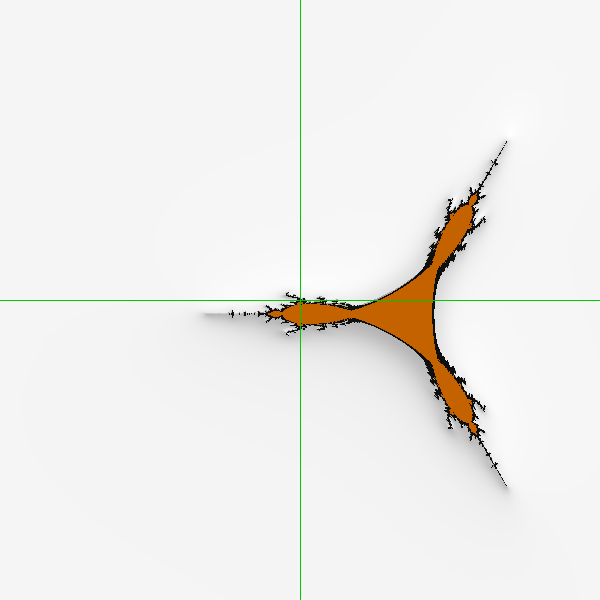}}
\fbox{\includegraphics[width=.18\textwidth, bb = 0 0 600 600]{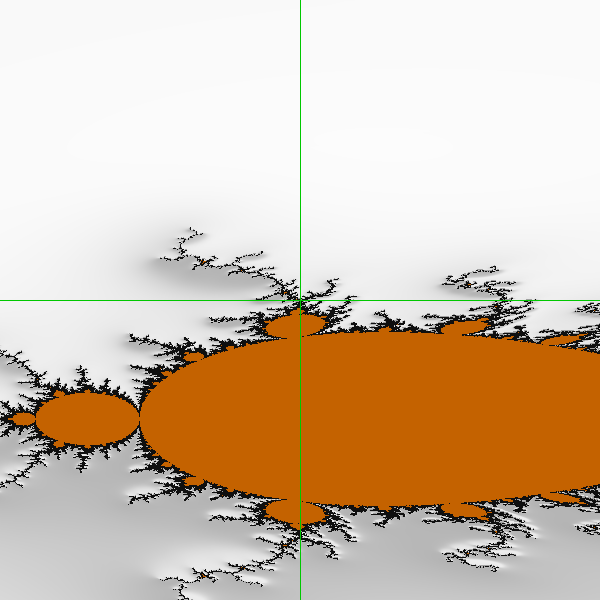}}
\fbox{\includegraphics[width=.18\textwidth, bb = 0 0 600 600]{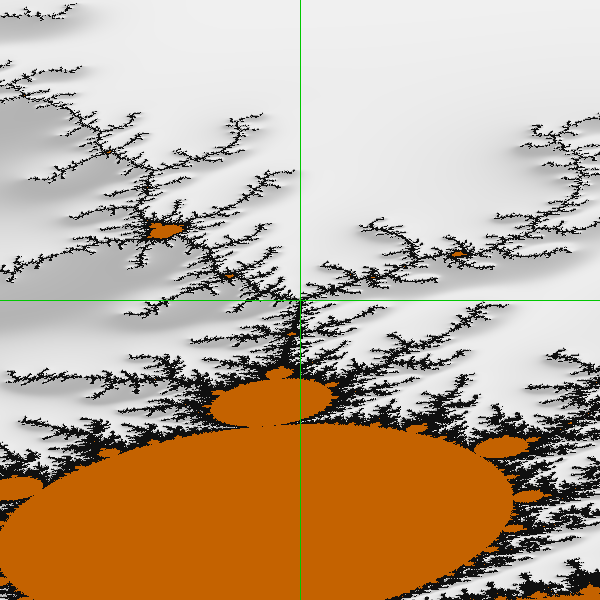}}
\fbox{\includegraphics[width=.18\textwidth, bb = 0 0 600 600]{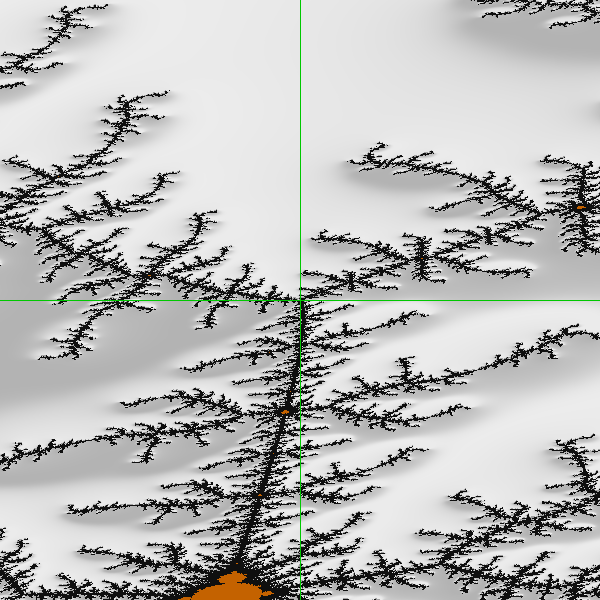}}
\fbox{\includegraphics[width=.18\textwidth, bb = 0 0 600 600]{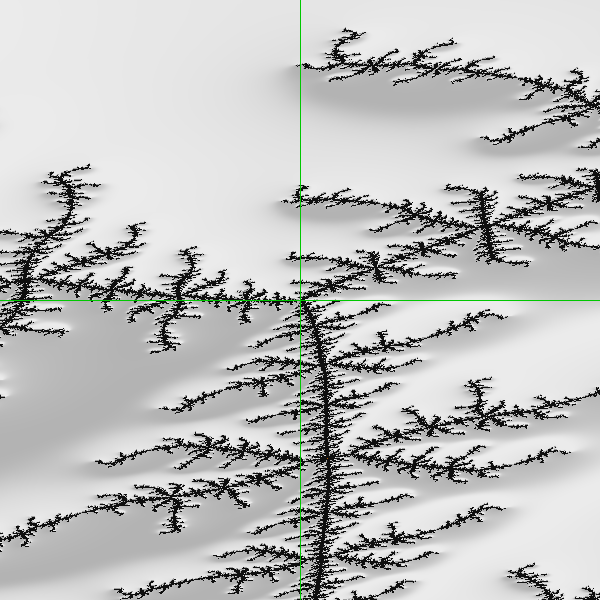}}
\\[.5em]
(JT2)
\end{center}
\caption{
(T) The tricorn $\T$.
(JT1) Center: $-1.2222454262925588+ 0.18411010266019595 i$, 
square width: from 6.0 to 0.005.
(JT2) Center: $-1.0672232757314006+0.13470887783195631 i$, 
square width: from 6.0 to 0.001. 
}
\label{fig_JT}
\end{figure}

%%%%%%----------table of pictures
\fboxsep=0pt
\fboxrule=.5pt
\begin{figure}[htbp]
\begin{center} 
\fbox{\includegraphics[width=.18\textwidth, bb = 0 0 600 600]{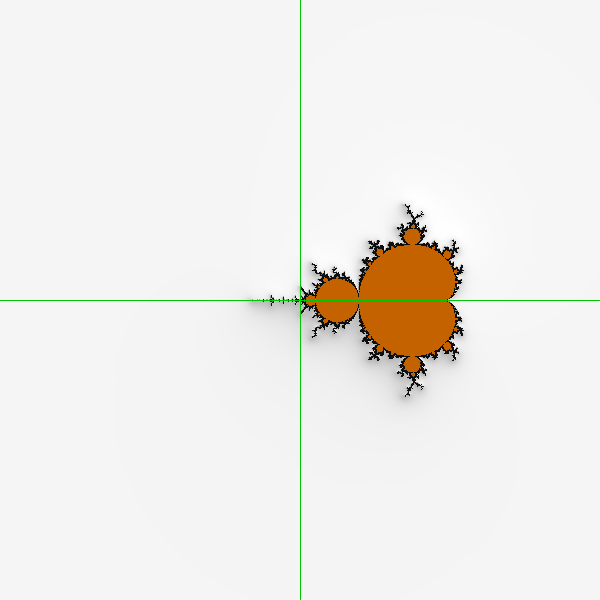}}
\fbox{\includegraphics[width=.18\textwidth, bb = 0 0 600 600]{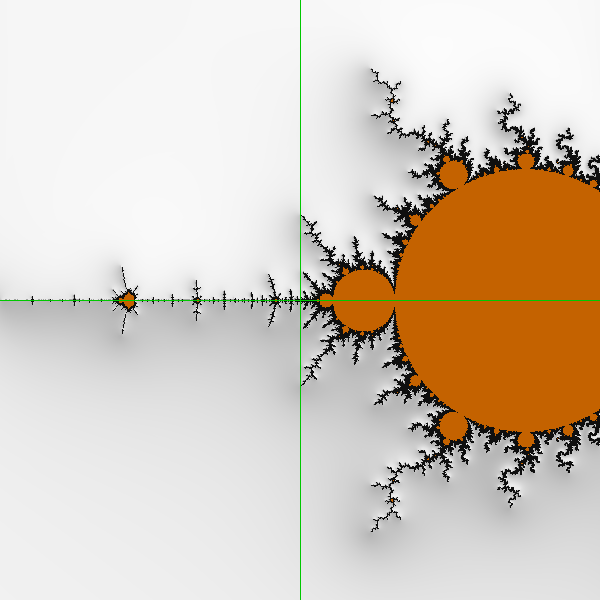}}
\fbox{\includegraphics[width=.18\textwidth, bb = 0 0 600 600]{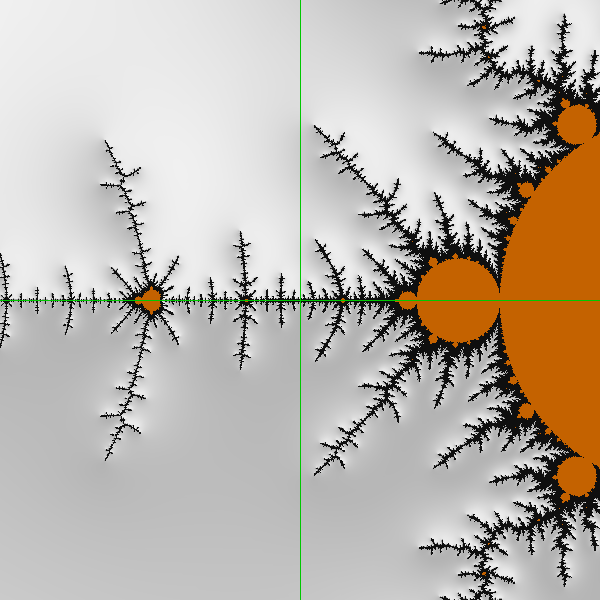}}
\fbox{\includegraphics[width=.18\textwidth, bb = 0 0 600 600]{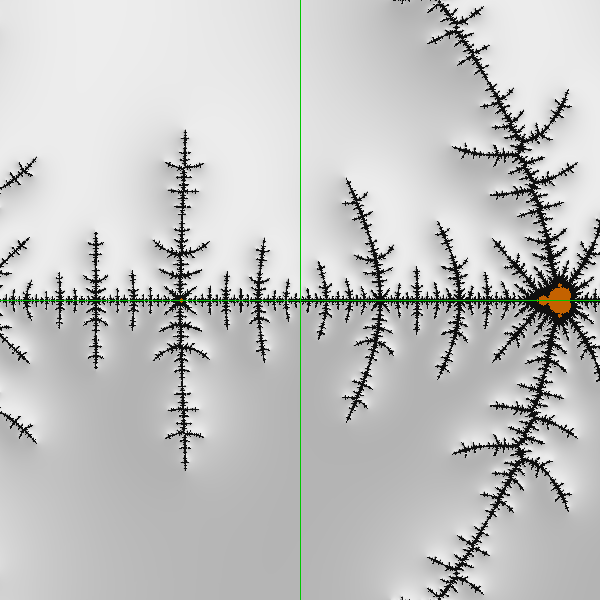}}
\fbox{\includegraphics[width=.18\textwidth, bb = 0 0 600 600]{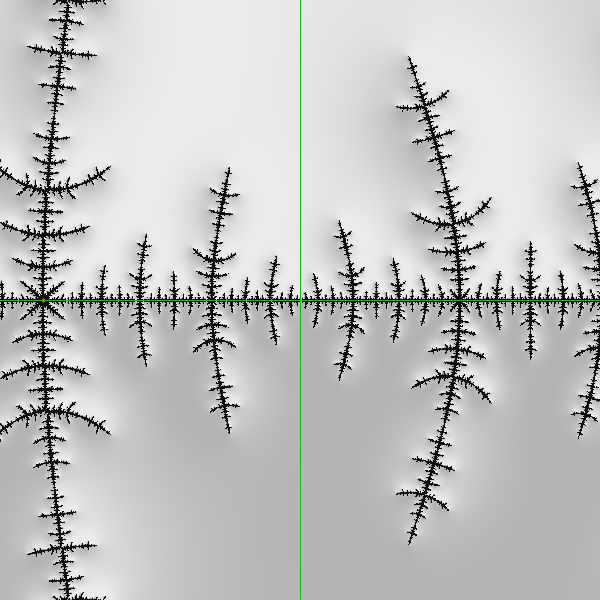}}
\\[.5em]
\fbox{\includegraphics[width=.18\textwidth, bb = 0 0 600 600]{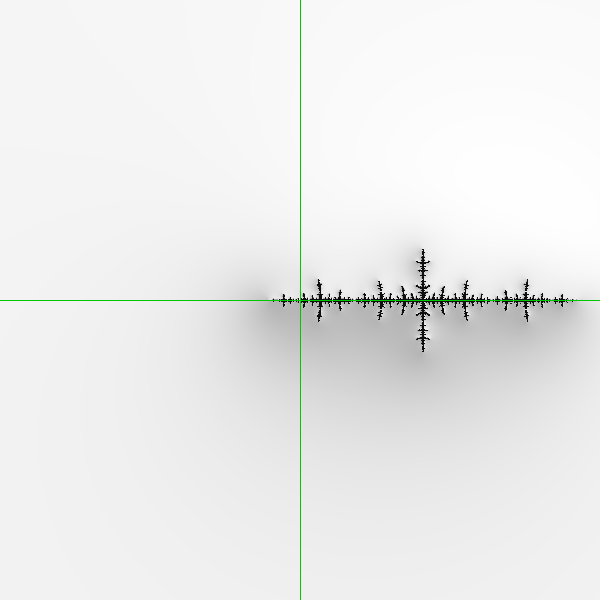}}
\fbox{\includegraphics[width=.18\textwidth, bb = 0 0 600 600]{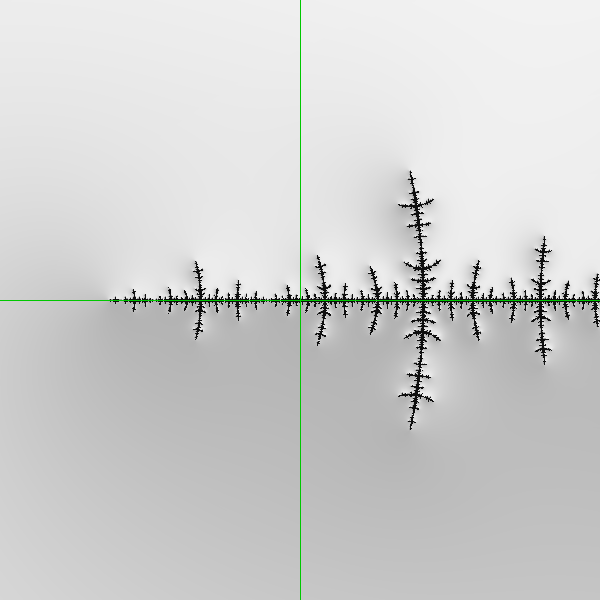}}
\fbox{\includegraphics[width=.18\textwidth, bb = 0 0 600 600]{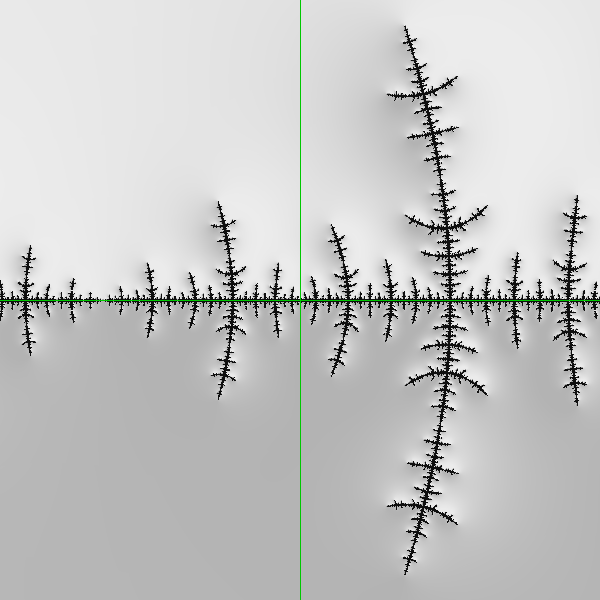}}
\fbox{\includegraphics[width=.18\textwidth, bb = 0 0 600 600]{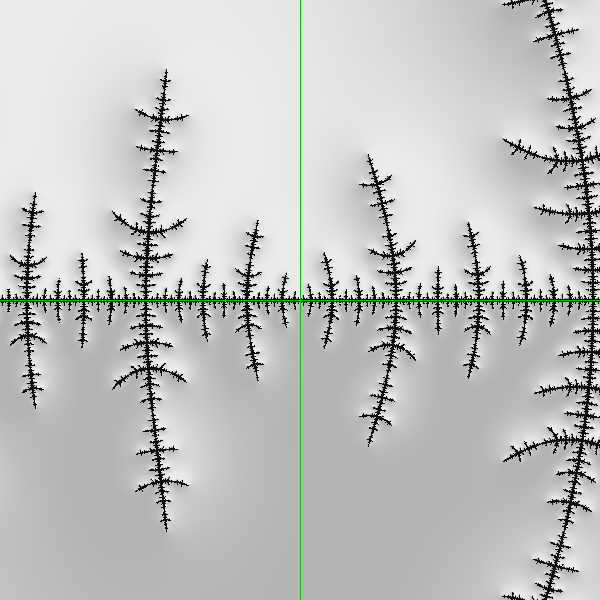}}
\fbox{\includegraphics[width=.18\textwidth, bb = 0 0 600 600]{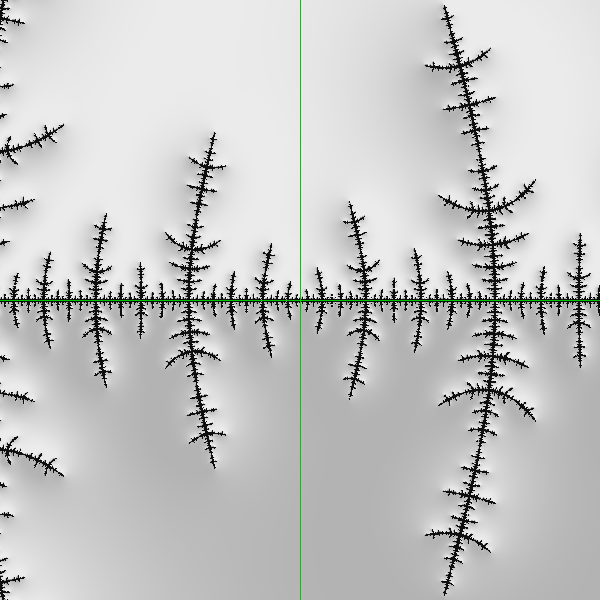}}
\\[.5em]
\fbox{\includegraphics[width=.18\textwidth, bb = 0 0 600 600]{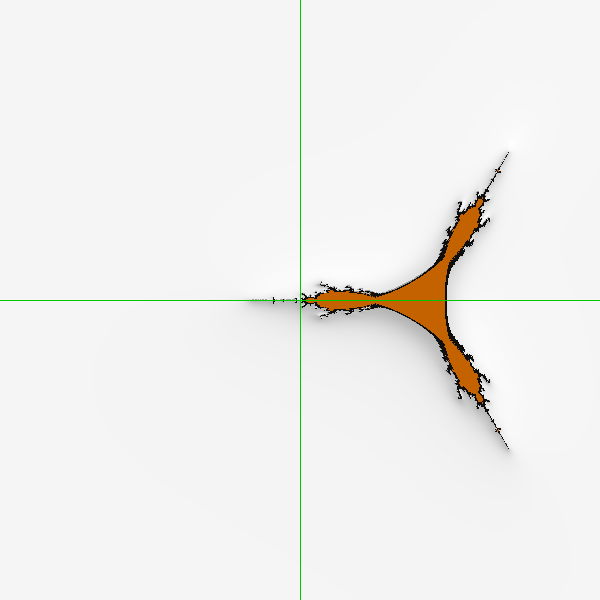}}
\fbox{\includegraphics[width=.18\textwidth, bb = 0 0 600 600]{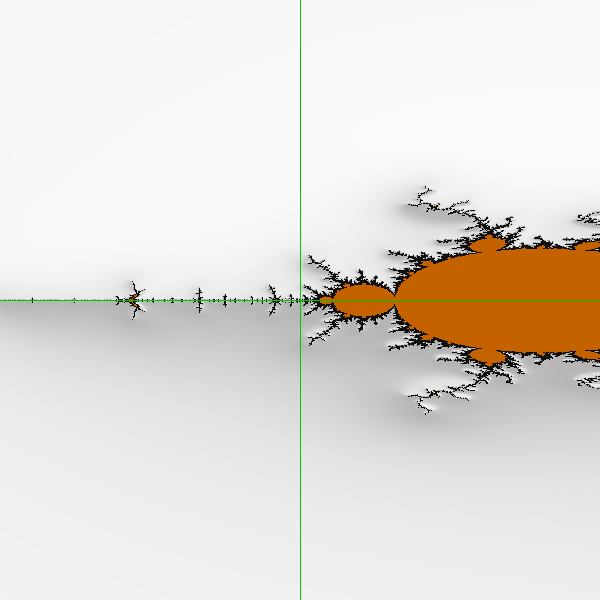}}
\fbox{\includegraphics[width=.18\textwidth, bb = 0 0 600 600]{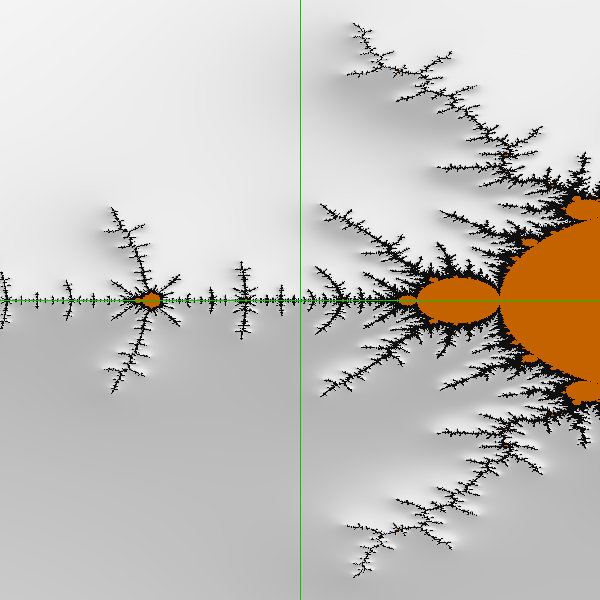}}
\fbox{\includegraphics[width=.18\textwidth, bb = 0 0 600 600]{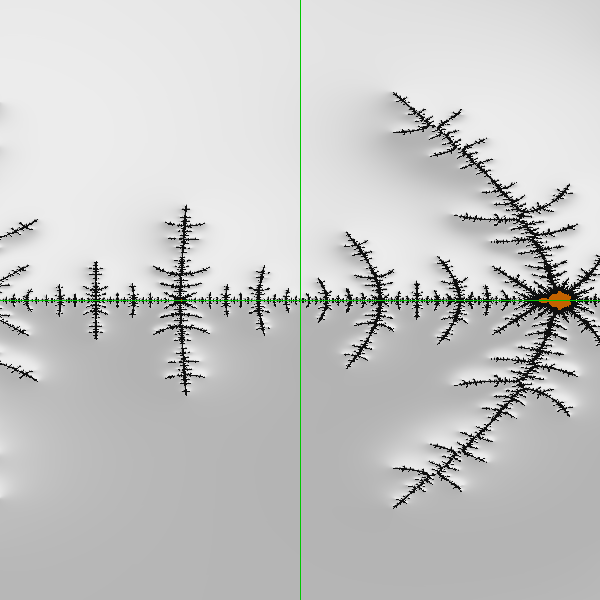}}
\fbox{\includegraphics[width=.18\textwidth, bb = 0 0 600 600]{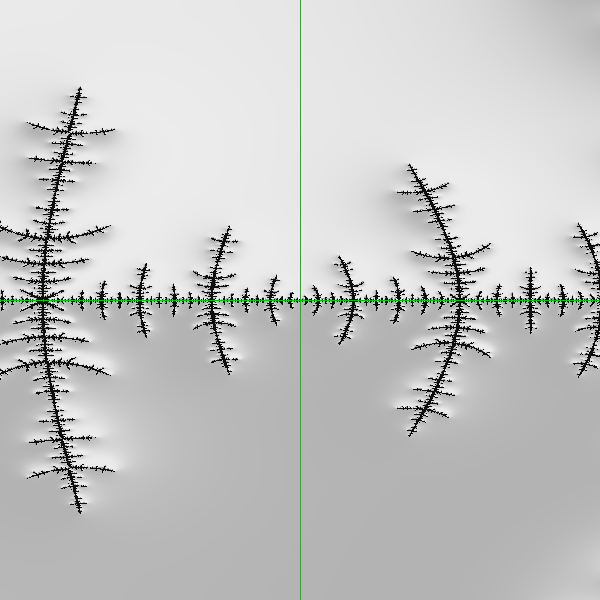}}
\\[.5em]
(MJT)
\end{center}
\caption{
(MJT) Center: $-1.4303576324513074$, square width: from 7.0 to 0.005.
}
\label{fig_MJT}
\end{figure}

\parag{Misiurewicz parameters.}
We say $c_0 \in \partial \T$ is a \text{\it Misiurewicz parameter} 
if the forward orbit of $c_0$ by $g_{c_0}$ eventually lands on 
a repelling periodic point. 
That is, there exist minimal $l \ge 1$ and $p \ge 1$ 
such that 
$g_{c_0}^l(c_0) = g_{c_0}^{l+p}(c_0)$ 
and $|Dg_{c_0}^p(g_{c_0}^l(c_0))| >1$.
(Since $g_{c_0}^p$ or $\overline{g_{c_0}^p}$ is 
holomorphic, the absolute value $|D{g_{c_0}^p}|$ is well-defined.)

Note that both $g_{c_0}^{2l}$ and $g_{c_0}^{2p}$ are holomorphic 
and that the relation $g_{c_0}^{2l}(c_0)= g_{c_0}^{2l+2p}(c_0)$
is satisfied. 
We let $a_0:=g_{c_0}^{2l}(c_0)$ and $\tilde{a}_0:=g_{c_0}(a_0)$,
which are repelling fixed points of $g_{c_0}^{2p}$ with the same 
multiplier $\lambda_0:=(g_{c_0}^{2p})'(a_0)$.

\parag{Remark.}
Unlike the Mandelbrot set, the Misiurewicz parameters are not dense 
in the boundary of $\T$. See \cite[Corollary 5.1]{IM}.

\parag{Similarity.}
Let us show an antiholomorphic counterpart of Tan's theorem (\thmref{thm_MJ}).
That is, the tricorn and the Julia sets are asymptotically similar 
at the Misiurewicz parameters {\it up to real linear transformation}.  
See Figures 2 and 3. In Figure 3, we take a real parameter and
compare the Mandelbrot set, the Julia set, and the tricorn. 

Let $c_0$ be a Misiurewicz parameter.
Then we have the following:

\begin{thm}[Similarity between $\T$ and $J^\ast$]\label{thm_TJ}
There exist an entire function $\phi$ on $\C$,
a real linear transformation $h:\C \to \C$,
and a sequence $\rho_k \to 0$
such that if we set $\cJ^\ast: = \phi^{-1}(J_{c_0}^\ast) \subset \C$,
then for any large constant $r>0$, we have
\begin{enumerate}[\rm (a)]
\item 
$\gauss{\rho_k^{-1}(J_{c_0}^\ast-c_0)}_r  \to [\cJ^\ast]_r$, and 
\item 
$\gauss{\rho_k^{-1}h(\T-c_0)}_r  \to [\cJ^\ast]_r$%~~(k \to \infty)$
\end{enumerate}
as $k \to \infty$ in the Hausdorff topology. 
\end{thm}

The proof is analogous to that of \thmref{thm_MJ}:
we start with showing an antiholomorphic 
version of \lemref{lem_key}. 
With the Misiurewicz parameter $c_0$,
integers $l$ and $p$, and repelling periodic point $a_0=g_{c_0}^{2l}(c_0)$ 
taken as above, we have the following:

\begin{lem}\label{lem_key_anti}
Suppose that $c_0 \in \T$ is a Misiurewicz parameter as above. 
For $k \in \N$, set $\rho_k:= 1/(g_{c_0}^{2l + 2kp})'(c_0)$.
Then we have the following.
\begin{enumerate}[\rm (1)]
\item
The function $\phi_k(w) = g_{c_0}^{2l + 2kp}(c_0 + \rho_k w)$ converges to a non-constant entire function 
$\phi:\C \to \C$ as $k \to \infty$ uniformly on any compact sets.
\item
There exists a real linear transformation $H(w)=Q w + Q' \overline{w}$ 
with $|Q| \neq |Q'|$
such that the function 
$$
\Phi_k(w) \dee 
g_{c_0 + H(\rho_k w)}^{2l + 2kp}(c_0 + H(\rho_k w))
$$
converges to the same function $\phi(w)$ as $k \to \infty$ 
uniformly on compact sets of $\C$.
\end{enumerate}
\end{lem}

Note that $\Phi_k(w)$ is a real analytic function of $w$
but it converges to an entire function $\phi(w)$.
We prove (1) and (2) separately.

\parag{Proof of (1).} 
Both $g_{c_0}^{2l}$ and $g_{c_0}^{2p}$ are holomorphic
and thus the usual derivatives 
$A_0 := (g_{c_0}^{2l})'(c_0) \neq 0$ 
and $\lam_0= (g_{c_0}^{2p})'(a_0)$ make sense. 
Since $a_0= g_{c_0}^{2l}(c_0)$ is a repelling fixed point
(that is, $|\lam_0|>1$) of $ g_{c_0}^{2p}$, the Poincar\'e function 
$$
\phi(w):= \lim_{n \to \infty} g_{c_0}^{2kp}
\paren{
a_0+\frac{w}{\lam_0^k}
}$$
exists and is an entire function,
where the convergence is uniform on compact sets. (See \thmref{thm_Poincare_func} in Appendix below.)
By the expansion $g_{c_0}^{2l}(c_0 + t)=a_0 +A_0t+o(t)~(t \to 0)$
we obtain
$$
\phi(w)= \lim_{n \to \infty} 
g_{c_0}^{2l+2kp}
\paren{c_0+\frac{w}{A_0 \lam_0^k}}.
$$
Hence we set $\rho_k:=(A_0\lam_0^k)^{-1}=1/(g_{c_0}^{2l + 2kp})'(c_0)$.
\QED

\medskip
To show (2), we need an extra lemma about stability 
and transversality of repelling periodic point $a_0$:

\begin{lem}[Stability and transversality]
\label{lem_stability_transversality}~
\begin{enumerate}[\rm (1)]
\item
Stability: There exists 
a real analytic function $c \mapsto a(c)$ defined near $c_0$
with $a(c_0)=a_0$ 
such that $a(c)$ is a repelling fixed point of 
$g_c^{2p}$, for which the multiplier 
$\lam(c):=(g_c^{2p})'(a(c))$
is also a real analytic function near $c_0$.
\item
Transversality: 
Let $b(c):=g_{c}^{2l}(c)$. 
Then there exist two complex numbers $B_0$ and $B_0'$ 
with $|B_0| \neq |B_0'|$ 
such that 
$$
b(c) -a(c) = B_0 (c-c_0)+ B_0' \overline{(c-c_0)}+o(|c-c_0|)
$$
as $c \to c_0$.
\end{enumerate}
\end{lem}
Note that
$B_0$ and $B_0'$ do not vanish simultaneously 
by the condition $|B_0| \neq |B_0'|$.

I learned the idea of using bi-quadratic maps 
to the antiholomorphic quadratic family 
from Hiroyuki Inou. 

\paragraph{Proof.}
Let $G_c(z):=g_c^2(z)=(z^2+\cbar)^2+c$, whose iterations are polynomial 
functions with variable $z$ with coefficients in $\Z[c, \overline{c}]$.
Then $a_0=g_{c_0}^{2l}(c_0)=G_{c_0}^l(c_0)$ is a repelling fixed point of $G_{c_0}^p$ 
with multiplier $(G_{c_0}^p)'(a_0)=\lam_0$. 

For (1), by a variant of the argument principle, 
such an $a(c)$ is explicitly given by
$$
a(c):= 
\frac{1}{2 \pi i}\int_{|z-a_0|=\e} 
z \cdot \frac{1-(G_c^p)'(z)}{z-G_c^p(z)} \dz
$$
where $z$ makes one turn anticlockwise around $a_0$
with a fixed, sufficiently small radius $\e$. 
Indeed, one can check that $G_c^p(a(c))=a(c)$ and
 real analyticity of $a(c)$ in $c$ comes from that of $G_c^p(z)$.
The multiplier $\lam(c):=(G_c^p)'(a(c))$ is real analytic as well,
and it satisfies $|\lam(c)|>1$ for $c$ sufficiently close to $c_0$,
since we have $|\lam(c_0)|>1$ by assumption.

To show (2) we employ a transversality result by van Strien
\cite[Main Theorem 1.1]{vS}, applied to the bi-quadratic family 
$$
\braces{F_{s,t}(z):=f_{t} \cc f_{s}(z)=(z^2+s)^2+t}_{s,\,t\, \in\, \C} \simeq \C^2.
$$
Since $g_{c}^2(z)=G_c(z)=F_{\cbar, c}(z)$,
our antiholomorphic family $\braces{g_c}_{c \in \C}$ can be regarded as 
a real analytic curve $\braces{(s,t)=(\cbar, c) \in \C^2 \st c \in \C}$
 in $\C^2$.
The critical points of $F_{s,t}$ are $\pm \sqrt{-s}$ and $0$,
but there is essentially one critical orbit when $(s,t)=(\cbar, c)$ 
since $F_{\cbar,c}(\pm \sqrt{-\cbar}) = g_c(0)=c$. 
One can also check that $F_{s,t}=A \cc F_{s',t'} \cc A^{-1}$ for some $A \in \Aff$
if and only if $(s',t')=(s,t),\,(\omega s, \omega^2 t)$ or $(\omega^2 s, \omega t)$.
This implies that the family $\{F_{s,t}\}$ is locally 
normalized near $(s,t)=(\overline{c_0}, c_0) \neq (0,0)$. 

When $(s,t)=(\overline{c_0}, c_0)$, we have
$$
a_0=G_{c_0}^{l+1}(\pm \sqrt{-\overline{c_0}})=
G_{c_0}^{l+p+1}(\pm \sqrt{-\overline{c_0}})
$$
and 
$$
\tilde{a}_0:=g_{c_0}(a_0)
=G_{c_0}^{l+1}(0)=
G_{c_0}^{l+p+1}(0),
$$
where both $a_0$ and $\tilde{a}_0$ are repelling fixed points of 
$G_{c_0}^p=F_{\overline{c_0}, c_0}^p$.
Then we have two analytic families 
$$
a(s,t):= 
\frac{1}{2 \pi i}\int_{|z-a_0|=\e} 
z \cdot \frac{1-(F_{s,t}^p)'(z)}{z-F_{s,t}^p(z)} \dz
$$
and
$$
\tilde{a}(s,t):= 
\frac{1}{2 \pi i}\int_{|z-\tilde{a}_0|=\e} 
z \cdot \frac{1-(F_{s,t}^p)'(z)}{z-F_{s,t}^p(z)} \dz
$$
of repelling fixed points of $F_{s,t}^p$, 
where $(s,t)$ are sufficiently close to $(\overline{c_0}, c_0)$ in $\C^2$.

By \cite[Main Theorem 1.1]{vS}, the map
$$
(s, t ) \mapsto 
\paren{
F_{s,t}^{l+1}(-\sqrt{-s})-a(s,t), F_{s,t}^{l+1}(\sqrt{-s})-a(s,t), 
F_{s,t}^{l+1}(0)-\tilde{a}(s,t)
} 
$$
is a local immersion near $(s,t)=(\overline{c_0}, c_0)$.
In other words, its Jacobian derivative at $(s,t)=(\overline{c_0}, c_0)$
has rank 2. 
Since the first and the second coordinates of the image always coincide, 
the derivative of the map 
$$
(s, t ) \mapsto 
(u,v):=
\paren{F_{s,t}^{l+1}(\sqrt{-s})-a(s,t), 
F_{s,t}^{l+1}(0)-\tilde{a}(s,t)}
$$
at $(s,t)=(\overline{c_0}, c_0)$ is rank 2 as well.
Hence if we write 
\begin{align*}
u&=B_0' (s-\overline{c_0}) + B_0 (t-c_0) + o\paren{\sqrt{(s-\overline{c_0})^2+(t-{c_0})^2}},
\quad\text{and} \\
v&=B_1' (s-\overline{c_0}) + B_1 (t-c_0) + o\paren{\sqrt{(s-\overline{c_0})^2+(t-{c_0})^2}},
\end{align*}
we have $B_0'B_1 - B_0B_1' \neq 0$. 
Now we let $(s,t)=(\overline{c}, c)$ with $c$ sufficiently close to $c_0$.
Then we have
$$
u=G_{c}^{l+1}(\sqrt{-\overline{c}})-a(\overline{c}, c)
=g_c^{2l}(c)-a(c)
=b(c)-a(c)
$$
and the expansion above implies 
$$
u=b(c)-a(c) = B_0' (\overline{c}-\overline{c_0}) + B_0 (c-c_0) +o(|c-c_0|)
$$
as in the statement. 
Hence it is enough to show $|B_0| \neq |B_0'|$.
Since we have 
$$
v=G_{c}^{l+1}(0)-\tilde{a}(\overline{c}, c)
=g_c^{2l+2}(0)-g_c(a(c))
=g_c(b(c))-g_c(a(c)),
$$
the expansion above implies 
$$
v=\overline{b(c)}^2-\overline{a(c)}^2 
= B_1' (\overline{c}-\overline{c_0}) + B_1 (c-c_0) +o(|c-c_0|).
$$
On the other hand, 
\begin{align*}
b(c)^2-a(c)^2 
&= (b(c)-a(c))(b(c)+a(c)) \\
&= \paren{B_0' (\overline{c}-\overline{c_0}) + B_0 (c-c_0) +o(|c-c_0|)}
\paren{2a_0+O(|c-c_0|)}\\
&=2 a_0 B_0' (\overline{c}-\overline{c_0}) + 2 a_0 B_0 (c-c_0) +o(|c-c_0|).
\end{align*}
Hence $B_1= \overline{2 a_0 B_0'}$ and  $B_1'= \overline{2 a_0 B_0}$
(where we have $a_0 \neq 0$, otherwise $a_0$ cannot be a repelling periodic point of $g_c$).
Since $B_0'B_1 - B_0B_1' \neq 0$, we conclude $|B_0'|^2-|B_0|^2 \neq 0$.
This completes the proof.
\QED

\parag{Remark.}
In the proof of Lemma 1 for the quadratic family,
we used the transversality result \cite[Lemma 1, p.333]{DH2},
which is proved in a purely algebraic way.
On the other hand, \cite[Main Theorem 1.1]{vS} used above 
assumes Thurston's rigidity theorem.
Recently Levin, Shen, and van Strien \cite{LSvS} 
showed a different version
of transversality (in the space of rational maps) 
by a rather elementary method.
It would be interesting to show \lemref{lem_stability_transversality}(2) 
in a purely algebraic way.

\medskip

Now let us go back to the proof of \lemref{lem_key_anti}.

\paragraph{Proof of \lemref{lem_key_anti}(2).}
Let us take a real linear transformation $H:\C \to \C$ 
of the form $W \mapsto H(W):=Q W+ Q'\overline{W}$,
with some complex constants $Q$ and $Q'$.
For $k \in \N$ and $w \in \C$ we define a real analytic function 
$\Phi_k(w)$ of the form 
$$
\Phi_k(w):= g_c^{2l+2kp}(c) = g_{c_0+H(\rho_k w)}^{2l+2kp}(c_0+H(\rho_k w)),
$$
where we let $c=c_0+H(\rho_k w)$ and $\rho_k=(A_0\lam_0^k)^{-1}$. 
We want to determine the constants $Q$ and $Q'$ such that
$\Phi_k(w)$ converges to $\phi(w)$ uniformly on compact sets of $\C$.

Let $\lam(c)=(g_c^{2p})'(a(c))$, as given in 
\lemref{lem_stability_transversality}(1).
Then the function $w \mapsto g_c^{2kp}(a(c)+w/\lam(c)^k)$ 
converges to the Poincar\'e function for $a(c)$ as $k \to \infty$.

On the other hand, 
when $c=c_0+H(\rho_k w)$ with $w$ in a compact set, we have 
\begin{align*}
\Phi_k(w)&=g_c^{2kp}(b(c))\\
&=g_c^{2kp}\paren{a(c)+B_0(c-c_0)+B_0' \overline{(c-c_0)}+o(|c-c_0|)}\\
&\sim g_c^{pk}\paren{a(c)+B_0 H(\rho_k w)+B_0' \overline{H(\rho_k w)}+o(\rho_k)}
\end{align*}
by \lemref{lem_stability_transversality}(2).
Hence it is enough to regard 
$B_0 H(\rho_k w)+B_0' \overline{H(\rho_k w)}$
as $w/\lam(c)^k$.
Since $\lam_0^k/\lam(c)^k \to 1$  as $k \to \infty$,
we obtain the condition
$$
B_0(Q\rho_k w +Q' \overline{\rho_k w})
+B_0' (\overline{Q\rho _k w} +\overline{Q'} \rho_k w)
=\frac{w}{\lam_0^k}.
$$
Since $\rho_k=1/(A_0\lam_0^k)$, we necessarily have
$$
Q=\frac{A_0\overline{B_0}}{|B_0|^2-|B_0'|^2}
\quad \text{and}\quad
Q'=-\frac{\overline{A_0}{B_0'}}{|B_0|^2-|B_0'|^2}.
$$
Note that $|B_0| \neq |B_0'|$ by \lemref{lem_stability_transversality}(2), 
and this implies $|B_0|^2-|B_0'|^2 \neq 0$ and $|Q| \neq |Q'|$.
Conversely, with these $Q$ and $Q'$ we have $\Phi_k(w) \to \phi(w)~(k \to \infty)$
uniformly on compact sets.
\QED  

\parag{Proof of \thmref{thm_TJ}.}
Let $h:=H^{-1}$, where $H$ is given in \lemref{lem_key_anti}.
Then the proof follows the same argument as that of \thmref{thm_MJ}.
However, an extra effort should be made when 
we choose a repelling periodic point $\phi(w_0')$ in the Julia set of $g_{c_0}$.
We should take a $w_0'$ (with the original condition $|w_0-w_0'| < \e /4$) 
such that $g_{c_0}^{2m}(\phi(w_0'))=\phi(w_0')$ 
for some $m \in \N$ and $\phi'(w_0) \neq 0$.
Then $w_0'$ is a simple zero of 
a holomorphic function $F(w):=g_{c_0}^{2m}(\phi(w))-\phi(w)$.
By \lemref{lem_key_anti}(2), this function is uniformly approximated by
a real analytic function (hence the Hurwitz theorem does not work)
$$
G_k(w):=g_{c_0+H(\rho_k w)}^{2m}(\Phi_k(w))-\Phi_k(w)
=g_c^{2m+2l+2kp}(c)-g_c^{2l+2kp}(c)
$$
near $w=w_0'$ with sufficiently large $k$, where $c=c_0+H(\rho_k w)$. 
Then one can show that $G_k$ maps a small round disk $D$ 
centered at $w_0'$
homeomorphically onto a topological disk containing $0$ 
when $k$ is sufficiently large;
 and that $w_k:=(G_k|_D)^{-1}(0)$ tends to $w_0'$
 as $k \to 0$. Then $c_k:=c_0+H(\rho_k w_k)$
 satisfies $g_{c_k}^{2m+2l+2kp}(c_k)=g_{c_k}^{2l+2kp}(c_k)$
 and this implies that $c_k \in \T$.
The remaining details are left to the reader.
\QED

\parag{Remark.}
\thmref{thm_TJ} can be easily generalized to 
the unicritical antiholomorphic family
$\braces{z \mapsto \zbar^d+c \st c \in \C}$ with $d \ge 2$. 
(In this case the tricorn becomes the {\it multicorn},
usually denoted by $\cM_d^\ast$.)

\appendix

\section*{Appendix. Existence of the Poincar\'e function} 
\label{sec_lin}
Here we give a proof of the existence of the Poincar\'e functions
associated with repelling periodic points.
This is originally shown by using a local linearization theorem by Koenigs.
See \cite[Cor.8.12]{Mi2}.
Our proof is based on the normal family argument and the univalent function theory (see \cite{Du} for example), which follows the idea of \cite[Lemma 4.7]{LM}.

\begin{thm}\label{thm_Poincare_func}
%Let $g$ be a holomorphic function 
%of the form $g(z) = \lam z +O(z^2)$ with $|\lam| >1$
%defined near $z = 0$. 
Let $g:\C \to \C$ be an entire function 
with $g(0) = 0,~g'(0) = \lam$, and $|\lam| >1$.
Then the sequence $\phi_n(w) = g^n(w/\lam^n)$ converges uniformly on compact sets in $\C$. 
Moreover, the limit function $\phi:\C \to \C$ satisfies $g \cc \phi(w) = \phi(\lam w)$ and $\phi'(0) = 1$.
\end{thm}

\parag{Proof.}
Since $g(z) = \lam z +O(z^2)$ near $z = 0$, there exists a disk $\Delta = \D(\delta) = \braces{z \in \C \st |z| < \delta}$ such that $g|\Delta$ is univalent and $\Delta \Subset g(\Delta)$. 
Hence we have a univalent branch $g_0^{-1}$ of $g$ that maps $\Delta$ into itself. 

First we show that $\phi_n$ is univalent on $\D(\delta/4)$: 
Since the map $\phi_n^{-1}:w \mapsto \lam^n g_0^{-n}(w)$ 
is well-defined on $\Delta = \D(\delta)$ and univalent, 
its image contains $\D(\delta/4)$ by the Koebe 1/4 theorem. 
Hence $\phi_n$ is univalent on $\D(\delta/4)$,
 and by the Koebe distortion theorem, 
 the family $\braces{\phi_n}_{n \ge 0}$ is 
 locally uniformly bounded on $\D(\delta/4)$ and thus equicontinuous. 

Next we show that $\phi_n$ has a limit on $\D(\delta/4)$:
Fix an arbitrarily large $r > 0$ and an integer $N$ such that $r < \delta |\lam|^N /4$.
By using the Koebe 1/4 theorem as above, 
the function $G_{N,k}(w):= \lam^N g^{k}(w/\lam^{N + k})~(k \in \N)$
satisfying $\phi_{N + k} =\phi_N \cc G_{N,k}$ 
is univalent on the disk $\D(\delta |\lam|^N/4)$. 
By the Koebe distortion theorem, 
there exists a constant $C>0$ independent of $N$ and $k$ such that 
for any $w \in \D(r)$ and sufficiently large $N$ 
we have $|G_{N,k}'(w)-1| \le C |w|/|\lam|^N$.
By integration we have $|G_{N,k}(w)-w| \le  C r^2/(2 |\lam|^N)$ on $\D(r)$. In particular, $G_{N,k} \to \id$ uniformly on $\D(\delta/4)$ as $N \to \infty$. 
Since the family $\braces{\phi_n}$ is equicontinuous 
on $\D(\delta/4)$, 
the relation $\phi_{N + k} =\phi_N \cc G_{N,k}$ implies that 
$\braces{\phi_n}_{n \ge 0}$ is Cauchy and has a unique limit $\phi$ 
on any compact sets in $\D(\delta/4)$. 

Let us check that the convergence extends to $\C$:
(We will not use the functional equation 
$g^n \cc \phi(w) = \phi(\lam^n w)$. Compare \cite[Cor.8.12]{Mi2}.)
Since $|\phi_{N + k}(w)-\phi_N(w)| = |\phi_N(G_{N,k}(w))-\phi_N(w)|$ 
and $|G_{N,k}(w)-w| = C r^2/(2 |\lam|^N)$ on $\D(r)$,
 it follows that the family $\braces{\phi_{N + k}}_{k \ge 0}$ (with fixed $N$) is uniformly bounded on $\D(r)$. 
Hence  $\braces{\phi_n}_{n \ge 0}$ is normal on any compact set in $\C$ and any sequential limit coincides with the local limit $\phi$ 
on $\D(\delta/4)$. 

The equation $g \cc \phi(w) = \phi(\lam w)$ and $\phi'(0) = 1$
are immediate from 
$g \cc \phi_{n}(w) = \phi_{n + 1}(\lam w)$ and $\phi_n'(0) = 1$.
\QED

\parag{Remark.} 
One can easily extend this proof to the case of meromorphic $g$
by using the spherical metric. 

\parag{Acknowledgement.}
The author would like to thank Hiroyuki Inou for suggesting 
the bi-quadratic family for the proof of 
\lemref{lem_stability_transversality}.
This work is partly supported by 
JSPS KAKENHI Grants Number 16K05193 and Number 19K03535.

{~}\\
Tomoki Kawahira\\
Department of Mathematics\\ 
Tokyo Institute of Technology\\
Tokyo 152-8551, Japan\\
kawahira@math.titech.ac.jp  \\

\noindent
Mathematical Science Team\\ 
RIKEN Center for Advanced Intelligence Project (AIP)\\
1-4-1 Nihonbashi, Chuo-ku\\ 
Tokyo 103-0027, Japan

\end{document}